\documentstyle[amssymb,12pt]{amsart}

\ifx\mathrm\undefined\let\mathrm\bf\fi
\ifx\mathbf\undefined\let\mathbf\bf\fi

\let\leq\leqslant \let\geq\geqslant
\let\tsize\textstyle \def\Sum{\sum\limits}

\newcommand{\half}{\frac12}
\newcommand{\Z}{{\Bbb Z}}

\newcommand{\C}{{\Bbb C}}
\newcommand{\Q}{{\Bbb Q}}
\newcommand{\Ref}[1]{{$($\ref{#1}$)$}}
\newcommand{\bean}{\begin{eqnarray}}
\newcommand{\eean}{\end{eqnarray}}
\newcommand{\be}{\begin{displaymath}}
\newcommand{\ee}{\end{displaymath}}
\newcommand{\bea}{\begin{eqnarray*}}
\newcommand{\eea}{\end{eqnarray*}}

\newcommand{\h}{{{\frak h\,}}}
\newcommand{\Id}{{\operatorname{Id}}}

\newcommand{\res}{{\operatorname{res}}}
\newcommand{\Ker}{{\operatorname{Ker}}}

\newcommand{\noi}{\noindent}
\newcommand{\vs}{\vspace{1.5\baselineskip}}
\newcommand{\Mu}{{\mathrm M}}
\newenvironment{proof}{\noindent{\it Proof\/}:}{$\;\Box$
\par\vs}
\newenvironment{definition}
{\noindent{\bf Definition\/}:}{\par\vs}
\newenvironment{example}
{\noindent{\bf Example\/}:}{\par\vs}
\newtheorem%
{thm}{Theorem}
\newtheorem%
{proposition}[thm]{Proposition}
\newtheorem%
{lemma}[thm]{Lemma}
\newtheorem%
{lemmadef}[thm]{Lemma-Definition}
\newtheorem%
{corollary}[thm]{Corollary}
\newtheorem%
{conjecture}[thm]{Conjecture}

\newcommand{\End}{{\operatorname{End}}}
\newcommand{\Hom}{{\operatorname{Hom}}}

\newcommand{\tr}{{\operatorname{tr}}}
\newcommand{\Sym}{{\operatorname{Sym}}}
\newcommand{\Fun}{{\operatorname{Fun}}}
\title[qKZB equations and Bethe ansatz]
{Solutions of the elliptic qKZB equations and Bethe ansatz I}
\thanks{${}^1$ Supported in part
by NSF grant DMS-9400841\\ \indent
$\hphantom{{}^1}$ Address after July 1, 1996: D-Math, ETH
Zentrum, CH-8092 Z\"urich, Switzerland}
\thanks{${}^2$ Supported in part by MAE--MICECO--CNRS Fellowship}
\thanks{${}^3$ Supported in part
by NSF grant DMS-9501290}
\author[G. Felder, V. Tarasov and A. Varchenko]
{G. Felder${}^{*,1}$, V. Tarasov${}^{**,2}$ and
A. Varchenko${}^{*,3}$}
\begin{document}
\maketitle
\vskip-.5\baselineskip
\centerline{${}^*${\it Department of Mathematics,
University of North Carolina at Chapel Hill,}}
\centerline{\it Chapel Hill, NC 27599-3250, USA}
\centerline{{\it E-mail addresses:} {\rm felder@@math.unc.edu,
varchenko@@math.unc.edu}}
\medskip\centerline{${}^{**}${\it Laboratoire de Physique Th\'eorique ENSLAPP,
\'Ecole Normale Sup\'erieure de Lyon,}}
\centerline{\it
46, All\'ee d'Italie, 69364 Lyon Cedex 07, France
}\centerline{{\it E-mail address:} {\rm {
vtarasov@@enslapp.ens-lyon.fr}}}
\bigskip
\medskip
\centerline{June 1996}
\bigskip
\medskip
\centerline{\it Dedicated to V. I. Arnold on his sixtieth birthday}
\medskip
\bigskip
\section{Introduction}
In this paper we give integral representations for solutions of
the system of elliptic quantum Knizhnik--Zamolodchikov--Bernard (qKZB)
difference equations in the case of $sl_2$. The qKZB equations
\cite{F} are
a quantum deformation of the KZB differential equations obeyed
by correlation functions of the Wess--Zumino--Witten model on
tori.
They have the form
\be
\Psi(z_1,\dots,z_j+p,\dots,z_n)=
K_j(z_1,\dots,z_n;\tau,\eta,p)\Psi(z_1,\dots,z_n).
\ee
The unknown function $\Psi$ takes values in a space of
vector valued functions of a complex variable $\lambda$, and
the $K_j$ are difference operators in $\lambda$. The parameters
of this system of equations are $\tau$ (the period of the elliptic curve),
$\eta$ (``Planck's constant''), $p$ (the step)
and $n$ ``highest weights''
$\Lambda_1,\dots,\Lambda_n\in\C$. The operators
$K_j$ are expressed in terms of $R$-matrices
of the elliptic quantum group $E_{\tau,\eta}(sl_2)$.

In the trigonometric limit $\tau\to i\infty$, the qKZB equations
reduce to the trigonometric qKZ equations \cite{FR}
obeyed by correlation
functions of statistical models and form factors of integrable
quantum field theories in 1+1 dimensions.

The KZB equations can be obtained in the semiclassical limit:
$\eta\to 0$, $p\to 0$, $p/\eta$ finite.

When the step $p$ of the qKZB equations goes to zero (with
the other parameters fixed) our construction gives common eigenfunctions
of the $n$ commuting operators $K_j(z_1,\dots,z_n;\tau,\eta,0)$
in the form of the Bethe
ansatz. These difference operators are closely related to the
transfer matrices of IRF models of statistical mechanics.

Our results follow from the main theme of this paper:
a geometric construction of tensor products of evaluation
Verma modules over the elliptic quantum group $E_{\tau,\eta}(sl_2)$.
In particular, we obtain some formulae given in
\cite{FV1} for the action of generators on these modules.

The results of this paper are thus parallel to the results on solutions
of the rational and trigonometric qKZB equations of \cite{V, TV2},
which are based on the representation theory of the Yangian
$Y(sl_2)$ and the affine quantum universal enveloping algebra
$U_q(\widehat{sl_2})$, respectively.

The paper is organized as follows: we begin by introducing the
notion of $R$-matrices and the qKZB equations in Section \ref{sqkzb}.
The geometric construction of $R$-matrices is given in
Section \ref{ssafs}. At the end of that section we show how to
obtain representations of $E_{\tau,\eta}(sl_2)$ in this way,
and give some explicit formulae for $R$-matrix elements.

The main applications of the geometric construction are contained
in Sections \ref{sba}--\ref{ssqkzb}. The Bethe ansatz eigenfunctions
of the commuting system of difference operators are given in
Section \ref{sba}, Theorem \ref{tba}. The proof of completeness
of Bloch eigenfunctions obtained by the Bethe ansatz
in the case of rational $\eta$ is given
for generic parameters in Section \ref{scbv}.
In Section \ref{caba}, we compare our results with the
algebraic Bethe ansatz of \cite{FV2}. In particular, we
show, using the results of Section
\ref{scbv} that the algebraic Bethe ansatz gives a basis of
eigenvectors in some cases.
Our results on
integral representation for
solutions of the qKZB equations are stated in Theorem \ref{irs0} and
Theorem \ref{irs}, Section \ref{ssqkzb}.
In Section \ref{sp}
we introduce the technique of iterated residues, the main
calculational tool of the paper, and complete the proofs of
our results.

Some of our results seem to be related to the results in
the recent preprint \cite{Ta}, in which integral solutions
of a system of difference equations associated to Sklyanin's algebra
are constructed.

In the next paper of this series, we will discuss
integration cycles, and compute the monodromy of the qKZB
equations. It turns out that the monodromy is described
in terms of $R$-matrices of
the elliptic quantum group $E_{p,\eta}(sl_2)$ and that
the elliptic quantum groups associated to the elliptic
curves $\C/(\Z+\tau\Z)$ and $\C/(\Z+p\Z)$ play a symmetric role
in the story.

\section{$R$-matrices, qKZB equations and
commuting difference operators}\label{sqkzb}

\subsection{$R$-matrices}
The qKZB equations are given in terms of $R$-matrices of
elliptic quantum groups. In the $sl_2$ case, these $R$-matrices
have the following properties. Let $\h=\C h$ be a one-dimensional
Lie algebra with generator $h$. For each $\Lambda\in\C$ consider the
$\h$-module $V_\Lambda=\oplus_{j=0}^\infty\C e_j$,
with $he_j=(\Lambda-2j)e_j$. For each pair $\Lambda_1$, $\Lambda_2$
of complex numbers we have a meromorphic function, called the $R$-matrix,
$R_{\Lambda_1,\Lambda_2}(z,\lambda)$ of two complex variables,
with values in $\End(V_{\Lambda_1}\otimes V_{\Lambda_2})$.

The main properties of the $R$-matrices are
\begin{enumerate}
\item[I.] Zero weight property: for any $\Lambda_i$, $z,\lambda$,
$[R_{\Lambda_1,\Lambda_2}(z,\lambda),h^{(1)}+h^{(2)}]=0$.
\item[II.] For any $\Lambda_1,\Lambda_2,\Lambda_3$,
the dynamical Yang--Baxter equation
\bea
&R_{\Lambda_1,\Lambda_2}(z,\lambda-2\eta h^{(3)})^{(12)}
R_{\Lambda_1,\Lambda_3}(z+w,\lambda)^{(13)}
R_{\Lambda_2,\Lambda_3}(w,\lambda-2\eta h^{(1)})^{(23)}
&
\\
&=
R_{\Lambda_2,\Lambda_3}(w,\lambda)^{(23)}
R_{\Lambda_1,\Lambda_3}(z+w,\lambda-2\eta h^{(2)})^{(13)}
R_{\Lambda_1,\Lambda_2}(z,\lambda)^{(12)},
\eea
holds in $\End(V_{\Lambda_1}\otimes V_{\Lambda_2}\otimes V_{\Lambda_3})$
for all $z,w,\lambda$.
\item[III.] For all $\Lambda_1$, $\Lambda_2$, $z,\lambda$,
$R_{\Lambda_1,\Lambda_2}(z,\lambda)^{(12)}
R_{\Lambda_2,\Lambda_1}(-z,\lambda)^{(21)}=\Id$. This property is
called ``unitarity''.
\end{enumerate}

We use the following notation: if $X\in\End(V_i)$
we denote by $X^{(i)}\in\End(V_1\otimes\dots\otimes V_n)$
the operator $\cdots\otimes\Id\otimes X\otimes\Id\otimes\cdots$, acting
non-trivially on the $i$th factor of a tensor product of vector spaces,
and
if $X=\sum X_k\otimes Y_k\in\End(V_i\otimes V_j)$ we set
$X^{(ij)}=\sum X_k^{(i)}Y_k^{(j)}$. If $X(\mu_1,\dots,\mu_n)$
is a function with values in $\End(V_1\otimes\dots\otimes V_n)$,
then $X(h^{(1)},\dots,h^{(n)})v=X(\mu_1,\dots,\mu_n)v$ if
$h^{(i)}v=\mu_iv$, for all $i=1,\dots,n$.

For each $\tau$ in the upper half plane and generic $\eta\in\C$
(``Planck's constant'')
a system of $R$-matrices $R_{\Lambda_1,\Lambda_2}(z,\lambda)$ obeying
I--III was constructed in \cite{FV1}. They are characterized
by an intertwining property with respect to the action of the
elliptic quantum group $E_{\tau,\eta}(sl_2)$ on tensor products
of evaluation Verma modules. In this paper,
we give an alternative geometric construction of these $R$-matrices,
see Section \ref{ssafs}.

\subsection{qKZB equations}
Fix the parameters $\tau,\eta$. Fix also
$n$ complex numbers $\Lambda_1,\dots,\Lambda_n$ and an additional
parameter $p\in \C$.
Let $V=V_{\Lambda_1}\otimes\cdots\otimes V_{\Lambda_n}$. The
kernel of $h^{(1)}+\dots+h^{(n)}$ on $V$ is called the zero-weight
space and is denoted $V[0]$. More generally, we write $V[\mu]$
for the eigenspace of $\sum h^{(i)}$ with eigenvalue $\mu$.
The qKZB equations
are difference equations for a function $\Psi(z_1,\dots,z_n,\lambda)$
of $n$ complex variables $z_1,\dots, z_n$ with values in the
space of meromorphic functions Fun$(V[0])$ of a complex variable $\lambda$
with values in $V[0]$.

The qKZB equations \cite{F} have the form
\bean\label{qKZB}
\Psi(z_1,\dots,z_j+p,\dots,z_n)=
R_{j,j-1}(z_j\!-\!z_{j-1}+p)\cdots
R_{j,1}(z_j\!-\!z_{1}+p)
\\
\Gamma_j
R_{j,n}(z_j\!-\!z_n)\cdots,R_{j,j+1}(z_j\!-\!z_{j+1})
\Psi(z_1,\dots,z_n)\notag
\eean
Here $R_{k,l}(z)$ is the operator of multiplication by
\be
\tsize R_{\Lambda_k,\Lambda_l}(z,\lambda-2\eta\Sum_{\tsize{j=1\atop j\neq k}}
^{l-1}h^{(j)})^{(k,l)}
\ee
acting
on the $k$th and $l$th factor of the tensor product, and
$\Gamma_j$ is the linear difference operator such that
$\Gamma_j\Psi(\lambda)=\Psi(\lambda-2\eta\mu)$ if $h^{(j)}\Psi=
\mu\Psi$.

The consistency of these equations follows from I--III.
In other words, the qKZB equations may be viewed as the
equation of horizontality for a flat discrete connection
on a trivial vector bundle with fiber Fun($V[0]$) over an
open subset of $\C^n$.

\subsection{Commuting difference operators}
A closely related set of difference equations is the eigenvalue
problem
\be
H_j(z)\psi=\epsilon_j\psi, \qquad j=1,\dots,n,\qquad \psi\in
{\mbox{Fun}}(V[0])
\ee
for the commuting difference operators
\begin{equation}\label{Go}
H_j(z)=R_{j,j-1}(z_j\!-\!z_{j-1})\cdots
R_{j,1}(z_j\!-\!z_{1})
\Gamma_j
R_{j,n}(z_j\!-\!z_n)\cdots,R_{j,j+1}(z_j\!-\!z_{j+1}).
\end{equation}
Here $z=(z_1,\dots,z_n)$ is a fixed generic point in $\C^n$ and $\psi$
is in Fun($V[0]$). The fact that the operators $H_j(z)$ commute with each
other follows from the flatness of the connection as $p\to 0$.

\subsection{Finite dimensional representations}
If $\Lambda$ is a nonnegative integer, $V_\Lambda$ contains the
subspace $SV_\Lambda=\oplus_{j=\Lambda+1}^\infty\C e_j$ with
the property that, for any $\Mu$,
$SV_\Lambda\otimes V_\Mu$ and $V_\Mu\otimes
SV_\Lambda$
are preserved by the $R$-matrices $R_{\Lambda,\Mu}(z,\lambda)$ and
$R_{\Mu,\Lambda}(z,\lambda)$, respectively, see \cite{FV1} and
Theorem \ref{tfd}. Let $L_\Lambda=V_\Lambda/SV_\Lambda$,
$\Lambda\in\Z_{\geq 0}$. Then, in particular, for any nonnegative integers
$\Lambda$ and $\Mu$,
$R_{\Lambda,\Mu}(z,\lambda)$ induces a map, also denoted by
$R_{\Lambda,\Mu}(z,\lambda)$, on the finite dimensional
space $L_\Lambda\otimes L_\Mu$.

The simplest nontrivial case is $\Lambda=\Mu=1$. Then
$R_{1,1}(z,\lambda)$
is defined on a four-dimensional vector space and coincides
with the {\em fundamental} $R$-matrix, the matrix of structure constants
of the elliptic quantum group $E_{\tau,\eta}(sl_2)$, see
Theorem \ref{tfrm}.

In any case, if $\Lambda_1,\dots,\Lambda_n$ are nonnegative integers,
we can consider the qKZB equations \Ref{qKZB} and the eigenvalue
problem \Ref{Go} on functions with values in the zero
weight space of $L_{\Lambda_1}\otimes
\cdots\otimes L_{\Lambda_n}$.

The results below obtained for the solutions with values in
$\otimes_j V_{\Lambda_j}$ immediately extend to this case:
let $\pi:\otimes_{j=1}^nV_{\Lambda_j}\to\otimes_{j=1}^{n}
L_{\Lambda_j}$ denote the canonical projection.

\begin{lemma}
Let $\Psi(z_1,\dots,z_n)$ be a solution of the qKZB equations
with values in $V[0]=V_{\Lambda_1}\otimes\cdots\otimes
V_{\Lambda_n}[0]$.
Then $\pi\circ\Psi(z_1,\dots,z_n)$ is a solution of the
qKZB equations with values in
$L[0]=L_{\Lambda_1}\otimes\cdots\otimes L_{\Lambda_n}[0]$.
Similarly if $\Psi\in\Fun(V[0])$ is a common eigenfunction
of the operators $H_j(z)$ then $\pi\circ\Psi$ is a common
eigenfunction of the induced operators $H_j(z)$ on $\Fun(L[0])$,
with the same eigenvalues.
\end{lemma}

In this lemma, the zero function should also be understood as an eigenfunction,
since, certainly, $\pi\circ\Psi$ can vanish.

\subsection{Remarks}\label{ssrem}
{\ }

\noi 1. If $\Lambda_1=\dots=\Lambda_n=\Lambda\in\Z_{\geq 0}$, then
the commuting operators $H_j(z)$ on $L[0]$ are special values
of the {\em transfer matrix} $T(w;z_1,\dots,z_n)$. This difference
operator is defined as follows. Let $P[\mu]\in\End(L_\Lambda)$ be
the projection onto the subspace of $L_\Lambda$ of weight $\mu$:
$P[\Lambda-2j]e_k=\delta_{jk}e_k$, $j=0,\dots,\Lambda$.
Define the partial trace
\be
\tr_{L_\Lambda[\mu]}=
(\tr\circ P[\mu])\otimes\Id:\End(\otimes_{j=0}^nL_\Lambda)\simeq
\End(L_\Lambda)\otimes\End(\otimes_{j=1}^nL_\Lambda)
\to\End(\otimes_{j=1}^nL_\Lambda)
\ee
Then the transfer matrix
$T(w)=T(w;z_1,\dots,z_n)\in\End(\Fun(L[0]))$ is defined by
\bea
T(w)f(\lambda)
&=&\tr_{L_\Lambda[\mu]}(R_{0,n}(w-z_n)\cdots
R_{0,1}(w-z_1))f(\lambda-2\eta\mu),
\\
R_{0,j}(w-z_j)&:=&R_{\Lambda,\Lambda}^{(0,j)}
(w-z_j,\lambda-2\eta\tsize\Sum_{k=1}^{j-1}
h^{(j)}).
\eea
These transfer matrices commute on $L[0]$
for different values of $w$.
It can be shown that
\be
H_j(z_1,\dots,z_n)=T(z_j;z_1,\dots,z_n).
\ee
In particular, the operators $H_j(z_1,\dots,z_n)$
commute with the transfer matrices.

\medskip
\noi 2. Below, a geometric construction of $R$-matrices is given.
This construction gives in particular a construction of
$R_{\Lambda,\Mu}$ on $L_\Lambda\otimes L_\Mu$
for all $\Lambda,\Mu\in\Z_{\geq 0}$. Alternatively, these $R$-matrices
can be computed starting from $R_{1,1}$ by the {\em fusion procedure},
see \cite{FV1}, Section 8.

\section{Modules over the elliptic quantum group
as function spaces}\label{ssafs}

In this section we realize the spaces dual to tensor products of
evaluation Verma modules over $E_{\tau,\eta}(sl_2)$ as spaces of functions.
The $R$-matrices are then constructed
geometrically.

Let us fix complex parameters $\tau$, $\eta$ with Im$(\tau)>0$,
and complex numbers $\Lambda_1,\dots, \Lambda_n$. We set
$a_i=\eta\Lambda_i$, $i=1,\dots,n$.

\subsection{A space of symmetric functions}\label{sssym}
We first introduce a space of functions with an action of the
symmetric group.
Recall that the Jacobi theta function
\begin{equation}\label{jacobi}
\theta(t)=-\sum_{j\in\Z}
e^{\pi i(j+\half)^2\tau+2\pi i(j+\half)(t+\half)},
\end{equation}
has multipliers $-1$ and $-\exp(-2\pi it-\pi i\tau)$ as $t\to t+1$
and $t\to t+\tau$, respectively. It is an odd entire function
whose zeros are simple and lie on the lattice $\Z+\tau \Z$. It
has the product formula
\be
\theta(t)\,=\,2e^{\pi i\tau/4}\sin(\pi t)
\prod_{j=1}^\infty(1-q^j)(1-q^je^{2\pi it})(1-q^je^{-2\pi it}),
\qquad q=e^{2\pi i\tau}.
\ee

\begin{definition} For complex numbers $a_1,\dots,a_n$,
$z_1,\dots,z_n$, $\lambda$
let $\tilde F^m_{a_1,\dots,a_n}(z_1,\dots,z_n,\lambda)$ be the space
of meromorphic functions $f(t_1,\dots,t_m)$ of
$m$ complex variables such that
\begin{enumerate}
\item[(i)] $\prod_{i<j}\theta(t_i-t_j+2\eta)
\prod_{i=1}^m\prod_{k=1}^n\theta(t_i-z_k-a_k)f$ is a
holomorphic function on $\C^m$.
\item[(ii)] $f$ is periodic with period 1 in each of its arguments
and
\be
f(\cdots,t_j+\tau,\cdots)
=e^{-2\pi i(\lambda+4\eta j-2\eta)}f(\cdots,t_j,\cdots),
\ee
for all $j=1,\dots,m$.
\end{enumerate}
\end{definition}
\noi There is an action of the symmetric group on this space of functions
that we introduce now.

\begin{lemma}\label{lsm}
The symmetric group $S_m$ acts on
$\tilde F^m_{a_1,\dots,a_n}(z_1,\dots,z_n,\lambda)$ so that
the transposition of $j$ and $j+1$ acts as
\be
s_jf(t_1,\dots,t_m)=f(t_1,\dots,t_{j+1},t_j,\dots,t_m)
\frac{\theta(t_{j}-t_{j+1}-2\eta)}
{\theta(t_{j}-t_{j+1}+2\eta)}\,.
\ee
\end{lemma}
\begin{proof}
Denote by $\phi(t_j-t_{j+1})$ the ratio of theta functions
in the definition of the action of $s_j$. The meromorphic
function $\phi$ is
1-periodic and obeys $\phi(t+\tau)=e^{8\pi i\eta}\phi(t)$.
Therefore, the action preserves the behavior of functions in $\tilde F^m$
under translations by the lattice. The position of the poles is also
preserved by the action as it is easy to check.
The relation $s_js_{j+1}s_j=s_{j+1}s_js_{j+1}$ holds automatically
and the relation $s_j^2=1$ follows from $\phi(t)\phi(-t)=1$.
These are the two relations defining the symmetric group.
\end{proof}

\begin{definition}
For any $m\in\Z_{>0}$,
let $F^m_{a_1,\dots,a_n}(z_1,\dots,z_n,\lambda)=\tilde F^m_{a_1,\dots,a_n}
(z_1,\dots,z_n,\lambda)^{S_m}$ be the space of $S_m$-invariant
functions. If $m=0$, we set $F^0_{a_1,\dots,a_n}(z_1,\dots,z_n,\lambda)=\C$.
We denote by Sym the symmetrization operator
Sym$\,=\,\sum_{s\in S_m}s:\tilde F^m\to F^m$.
Also, we set
\be
F_{a_1,\dots,a_n}(z_1,\dots,z_n,\lambda)
=\oplus_{m=0}^\infty F^m_{a_1,\dots,a_n}(z_1,\dots,z_n,\lambda),
\ee
and define an $\h$-module structure on
$F_{a_1,\dots,a_{n}}(z_1,\dots,z_{n},\lambda))$ by letting $h$ act
by
\be
h|_{F^{m}_{a_1,\dots,a_{n}}(z_1,\dots,z_{n},\lambda)}=(\tsize\Sum_{i=1}^n
\Lambda_i-2m)\Id,\qquad a_i=\eta\Lambda_i.
\ee
\end{definition}

Some technical results on this space of functions are
proved in Section \ref{sp}. In particular, it is shown
in Lemma \ref{ldim} that
$F^{m}_{a_1,\dots,a_{n}}(z_1,\dots,z_{n},\lambda)$ is
a finite dimensional vector space of dimension
$\left(\begin{matrix}{n+m-1}\\ {m}\end{matrix}\right)$.

Clearly, $F^m_{a_{\sigma(1)},\dots,a_{\sigma(n)}}(z_{\sigma(1)},\dots,
z_{\sigma(n)})=F^{m}_{a_1,\dots,a_{n}}(z_1,\dots,z_{n},\lambda)$
for any permutation $\sigma\in S_n$.

\begin{example}
Let $n=1$. Then
$F^m_a(z,\lambda)$ is a one-dimensional space spanned by
\begin{equation}\label{eb}
\omega_m(t_1,\dots,t_m,\lambda;z)=
\prod_{i<j}\frac{\theta(t_i-t_j)}
{\theta(t_i-t_j+2\eta)}\prod_{j=1}^{m}\frac
{\theta(\lambda+2\eta m+t_j-z-a)}
{\theta(t_j-z-a)},
\end{equation}
see Lemma \ref{ln1}, Section \ref{sp}.
\end{example}

\subsection{Tensor products}\label{sstp}
\begin{proposition}\label{p32}
Let $n=n'+n''$, $m=m'+m''$ be nonnegative integers
and $a_1,\dots,a_n$, $z_1,\dots,z_n$ be complex numbers.
The formula
\be
k(t_1,\dots,t_{m})=
\frac1{m'!m''!}\Sym \biggl( f(t_{1},\dots,t_{m'})g(t_{m'+1},\dots,t_{m})
\prod_{
\begin{matrix} \scriptstyle{m'<j\leq m}\\
\scriptstyle{1\leq l\leq n'}\end{matrix}}
\frac{\theta(t_j-z_l+a_l)}{\theta(t_j-z_l-a_l)}
\biggr)\ee
correctly defines a linear map $\Phi:f\otimes g\mapsto k=\Phi(f\otimes g)$,
\bea
\oplus_{m'=0}^{m}
F^{m'}_{a_1,\dots,a_{n'}}(z_1,\dots,z_{n'},\lambda)
\otimes
F^{m''}_{a_{n'+1},\dots,a_{n}}(z_{n'+1},\dots,z_{n},\lambda-2\nu)
\\
\to
F^{m}_{a_1,\dots,a_{n}}(z_1,\dots,z_{n},\lambda),
\eea
where
$
\nu=a_{1}+\cdots +a_{n'}-2\eta m'.
$
For generic values of the parameters $z_j$, $\lambda$, the map $\Phi$
is an isomorphism. Moreover, $\Phi$ is associative in the
sense that, for any three functions $f,g,h$,
$\Phi(\Phi(f\otimes g)\otimes h)=
\Phi(f\otimes\Phi(g\otimes h))$, whenever defined.
\end{proposition}

The only claim that does not follow immediately
{}from the definitions is the claim that $\Phi$ is an isomorphism
for generic $z_i$'s and $\lambda$. The proof of this is deferred
to Section \ref{sp}.

By iterating this construction, we get for all $n\geq1$
a linear map $\Phi_n$,
defined recursively by $\Phi_1=\Id$,
$\Phi_{n}=\Phi(\Phi_{n-1}\otimes\Id)$, from
\be
\oplus_{m_1+\cdots+m_n=m}\otimes_{i=1}^n
F_{a_i}^{m_i}(z_i,\lambda-2\eta(\mu_{1}+\cdots+\mu_{i-1}))
\ee
to $F^m_{a_1,\dots,a_n}(z_1,\dots,z_n,\lambda)$, with
$\mu_j=a_j/\eta-2m_j$, $j=1, \dots, n$.
Let $V_\Lambda^*=\oplus_{j=0}^\infty\C e_j^*$ be the restricted dual of
the
module $V_\Lambda=\oplus_{j=0}^\infty \C e_j$. It is spanned by the
basis $(e_j^*)$ dual to the basis $(e_j)$.
We let $\h$ act on
$V_\Lambda^*$ by $he_j^*=(\Lambda-2j)e_j^*$.
Then the map that sends $e^*_j$ to
$\omega_j$ (see \Ref{eb}) defines an
isomorphism of $\h$-modules
\be
\omega(z,\lambda):V_\Lambda^*\to F_a(z,\lambda), \qquad a=\eta\Lambda.
\ee
By composing this with the maps $\Phi$ of Proposition \ref{p32}, we obtain
homomorphisms (of $\h$-modules)
\be
\omega(z_1,\dots,z_n,\lambda):
V_{\Lambda_1}^*\otimes\cdots\otimes
V_{\Lambda_n}^*\to F_{a_1,\dots,a_n}(z_1,\dots,z_n,\lambda)
\ee
which are isomorphisms for generic values of $z_1,\dots,z_n,\lambda$.
The restriction of the map $\omega(z_1,\dots,z_n,\lambda)$
to $\C e_{m_1}^*\otimes\cdots\otimes e_{m_n}^*$
is
\be
\Phi_n(\omega(z_1,\lambda)\otimes
\omega(z_{2},\lambda-2\eta\mu_1)\otimes
\cdots\otimes\omega(z_n,\lambda-2\eta(\mu_1+\cdots+\mu_{n-1}))),
\ee
where $\mu_j=\Lambda_j-2m_j$, $j=1,\dots,n$.
For example, if $n=2$, then $\omega(z_1,z_2,\lambda)$ sends
$e^*_j\otimes e^*_k$ to
\be
\frac1{j!k!}\Sym\biggl(
\omega_j(t_1,\dots,t_{j},\!\lambda;z_1)
\omega_k(t_{j+1},\dots,t_{j+k},\!\lambda-2a_1+4\eta j;z_2)\!
\prod_{i=j+1}^{j+k}\!\frac{\theta(t_i-z_1+
a_1)}{\theta(t_i-z_1-a_1)}\biggr),
\ee
where $\{\omega_j(t_1,\dots,t_j,\lambda;z)\}$ is the basis \Ref{eb}
of $F_a(z,\lambda)$.

More generally, we have an explicit formula for
the image of $e_{m_1}^*\otimes\cdots\otimes e_{m_n}^*$, which
we discuss next.

\subsection{A basis of $F_{a_1,\dots,a_n}(z_1,\dots,z_n)$}
The space $V_\Lambda$ comes with a basis $e_j$. Thus we have the
natural basis $e_{m_1}^*\otimes\cdots\otimes e_{m_n}^*$ of the
tensor product of $V^*_{\Lambda_i}$ in terms of the dual bases
of the factors.
The map $\omega(z_1,\dots,z_n,\lambda)$ maps, for generic $z_i$,
this basis to a basis of $F_{a_1,\dots, a_n}(z_1,\dots,z_n,\lambda)$,
which is an essential part of our formulae for Bethe
ansatz eigenvectors and integral representations for solutions
of the qKZB equations.

We give here an explicit formula for the basis vectors.
\begin{proposition}\label{pomegam}
Let $m\in\Z_{\geq 0}$, $\Lambda=(\Lambda_1,\dots,\Lambda_n)\in\C^n$,
and let $z=(z_1,\dots,z_n)\in\C^n$ be
generic. Set $a_i=\eta\Lambda_i$.
Let
\be
u(t_1,\dots,t_m)=\prod_{i<j}
\frac{\theta(t_i-t_j+2\eta)}
{\theta(t_i-t_j)}
\ee
Then, for
generic $\lambda\in\C$,
the functions
\be
\omega_{m_1,\dots,m_n}(t_1,\dots,t_m,\lambda;z)=\omega(z,\lambda)
\,e_{m_1}^*\otimes\cdots\otimes e_{m_n}^*
\ee
labeled by $m_1,\dots,m_n\in\Z$ with $\sum_km_k=m$
form a basis of $F^m_a(z,\lambda)$ and are given by the
explicit formula
\begin{gather*}
\omega_{m_1,\dots,m_n}(t_1,\dots,t_m,\lambda;z)=
u(t_1,\dots,t_m)^{-1}
\sum_{I_1,\dots,I_n}
\prod_{l=1}^n
\prod_{i\in I_l}
\prod_{k=1}^{l-1}
\frac
{\theta(t_i-z_k+a_k)}
{\theta(t_i-z_k-a_k)}
\\
\times\prod_{k<l}
\prod_{i\in I_k,j\in I_l}
\frac
{\theta(t_i-t_j+2\eta)}
{\theta(t_i-t_j)}
\prod_{k=1}^{n}
\prod_{j\in I_k}
\frac
{\theta(\lambda\!+\!t_j\!-\!z_k\!-\!a_k\!+\!2\eta
m_k\!-\!2\eta\sum_{l=1}^{k-1}(\Lambda_l\!-\!2m_l))}
{\theta(t_j-z_k-a_k)}\, .
\end{gather*}
The summation is over all $n$-tuples $I_1,\dots,I_n$
of disjoint subsets of $\{1,\dots,m\}$ such that
$I_j$ has $m_j$ elements, $1\leq j\leq n$.
\end{proposition}

\begin{proof} If $n=1$ this formula is the same as the one given
above, eq.\ \Ref{eb}.
The function $u$ is designed to intertwine between the twisted
symmetrization Sym defined above and the ordinary symmetrization
Sym${}_0:f\mapsto\sum_{\sigma\in S_m}f(\sigma t)$. Indeed, for any function
$f$ of $t_1,\dots,t_m$, we have $\Sym_0 (uf)
=u\,\Sym(f)$.
Let $I^0_1=\{1,\dots,m_1\}$, $I^0_2=\{m_1+1,\dots,m_1+m_2\}$ and
so on. Also, if $I=\{i_1<\cdots<i_k\}$ set
$t_I=(t_{i_1},\dots,t_{i_k})$.
Then, by definition, $\omega_{m_1,\dots, m_n}$ is given by
\be
\frac1{m_1!\cdots m_n!}
\Sym \left(
\prod_{k=1}^n
\omega_{m_k}(t_{I^0_{k}},\lambda-2\eta{\tsize\Sum_{l=1}^{k-1}}
(\Lambda_l-2m_l);z_k)\prod_{l=1}^{k-1}\prod_{j\in I^0_l}
\frac{\theta(t_j-z_l+a_l)}
{\theta(t_j-z_l-a_l)}\right).
\ee
The one-point functions $\omega_{m_k}$ are given in eq.\ \Ref{eb}.
Therefore, $(\prod m_j!)u\,\omega_{m_1,\dots,m_n}$ can be computed by ordinary
symmetrization of the product of $u$ with the expression in the
brackets. It is then straightforward to compute this product, which
is given by the term with $I_j=I_j^0$ in the sum above. The
symmetrization gives all terms in the sum, each with multiplicity
$\prod m_j!$. Thus the factorials are canceled and we get the
claimed formula.
\end{proof}

\subsection{$R$-matrices}
Let $a=\eta\Lambda$ and $b=\eta \Mu$ be complex numbers.
Since $F_{ab}(z,w,\lambda)=F_{ba}(w,z,\lambda)$ by the symmetry of the
definition, we obtain a family of isomorphisms between $V^*_\Lambda \otimes
V^*_\Mu$ and $V^*_\Mu\otimes V^*_\Lambda $. The composition of this family with
the flip $P:V^*_\Mu\otimes V^*_\Lambda \to V^*_\Lambda \otimes V^*_\Mu$,
$Pv\otimes w=w\otimes v$ gives
a family of automorphisms of $V^*_\Lambda \otimes V^*_\Mu$:

\vs\begin{definition}
Let $z,w,\lambda$ be such that $\omega(z,w,\lambda)
:V^*_\Lambda \otimes V^*_\Mu\to F_{ab}(z,w,\lambda)$ is invertible.
The {\em $R$-matrix}
$R_{\Lambda,\Mu}(z,w,\lambda)\in\End_\h(V_\Lambda \otimes V_\Mu)$
is the dual map to the composition $R^*_{\Lambda,\Mu}(z,w,\lambda)$:
\be
V^*_\Lambda\otimes V^*_\Mu
\buildrel{P}\over{\longrightarrow}
V^*_\Mu\otimes V^*_\Lambda
\buildrel{\omega(w,z,\lambda)}\over{\longrightarrow}
F_{ab}(z,w,\lambda)
\buildrel{\omega(z,w,\lambda)^{-1}}\over{\longrightarrow}
V^*_\Lambda \otimes V^*_\Mu,
\ee
where we identify canonically $V_\Lambda ^*\otimes V_\Mu^*$ with
$(V_\Lambda \otimes V_\Mu)^*$.
\end{definition}

Alternatively, the $R$-matrix $R_{\Lambda,\Mu}(z,w,\lambda)$
can be thought of as a transition matrix
expressing the basis $\tilde\omega_{ij}=
\omega(w,z,\lambda)e_j^*\otimes e_i^*$ in terms of the basis
$\omega_{ij}=\omega(z,w,\lambda)e_i^*\otimes e_j^*$ of $F_{ab}(z,w,\lambda)$:
if $R_{\Lambda, \Mu}(z,w,\lambda)e_i\otimes e_j=\sum_{kl}
R_{ij}^{kl}e_k\otimes e_l$, then
\be
\tilde\omega_{kl}=\sum_{ij}R_{ij}^{kl}\omega_{ij}.
\ee

\medskip
\noi{\bf Remark.} Note that by the symmetry
under reversal of all
signs,
\be
R_{\Lambda,\Mu}(z,\lambda)=R_{-\Lambda,-\Mu}(-z,-\lambda),
\ee
if we identify $V_\Lambda$ with $V_{-\Lambda}$ via the basis
($e_j$).

\begin{lemma}\label{lpp}
\begin{enumerate}
\item[(i)] $R_{\Lambda, \Mu}(z,w,\lambda)$ is a meromorphic function
of $\Lambda, \Mu,z,w,\lambda$.
\item[(ii)] If $\Lambda$ is generic, then
$R_{\Lambda,\Lambda}(z,w,\lambda)$ is regular at $z=w$ and
$\lim_{z\to w}R_{\Lambda,\Lambda}(z,w,\lambda)=P$, the flip $u\otimes v\mapsto
v\otimes u$.
\item[(iii)]
$R_{\Lambda, \Mu}(z,w,\lambda)$ depends only on the difference
$z-w$.
\end{enumerate}
\end{lemma}

\begin{proof} (i)
Consider the bases $\{\omega_{ij}\}$, $\{\tilde \omega_{ij}\}$,
$i+j=m$,
of $F^m_{ab}(z,w,\lambda)$. These functions are meromorphic in
all variables. Since they are linearly independent there
are $m+1$ values $t^0,\dots,t^m$
of $t\in\C^{m}$ so that $(\omega_{i,m-i}(t^j))_{0\leq i,j\leq m}$
is an invertible matrix $B$. The $R$-matrix is then
$B^{-1}\tilde B$ where $\tilde B$ is the matrix with entries
$\tilde\omega_{i,m-i}(t^j)$. Therefore it is meromorphic.

\noi (ii) The proof of this is deferred to the end of
Section \ref{sp}.

\noi(iii)
We have, for all $c\in\C$,
an isomorphism
$U_c:F_{ab}(z,w,\lambda)\to F_{ab}(z+c,w+c,\lambda)$, given by
the translation of the arguments $t_i$ by $c$. It is easy
to see that
$\omega(z+c,w+c,\lambda)U_c=\omega(z,w,\lambda)$, so
$R_{\Lambda, \Mu}(z+c,w+c,\lambda)=R_{\Lambda, \Mu}(z,w,\lambda)$.
\end{proof}

Accordingly, we write $R_{\Lambda, \Mu}(z-w,\lambda)$ instead of
$R_{\Lambda, \Mu}(z,w,\lambda)$ in what follows.

By definition the $R$-matrix is determined by the
relation $\omega(z_1,z_2,\lambda)R_{\Lambda, \Mu}^*(z_1-z_2,\lambda)
=\omega(z_2,z_1,\lambda)P$ in $\Hom(V_\Lambda ^*\otimes V_\Mu^*,
F_{ab}(z_1,z_2,\lambda))$. More generally, by the
associativity of $\Phi$, we have:

\begin{lemma}\label{lsym}
For any $\Lambda_1,\dots,\Lambda_n$, the identity
\bea
\omega(\dots,z_{j},z_{j+1},\dots,\lambda)\!\!&\!\!R^*_{\Lambda_j\Lambda_{j+1}}
(z_j\!-\!z_{j+1},\lambda-2\eta\sum_{l<j}h^{(l)})^{(j,j+1)}
\\
&=\omega(\dots,z_{j+1},z_{j},\dots,\lambda)P^{(j,j+1)}&
\eea
holds in $\Hom(V^*_{\Lambda_1}\otimes\cdots\otimes V^*_{\Lambda_n},
F_{a_1,\dots,a_n}(z_1,\dots,z_n,\lambda))$, $a_i=\eta\Lambda_i$.
\end{lemma}

A corollary of this lemma is the
dynamical Yang--Baxter equation for the $R$-matrices:

\begin{thm}
The matrices $R_{\Lambda,\Mu}(z,\lambda)$ obey
I--III of Section \ref{sqkzb}.
\end{thm}
\begin{proof} Fix $z_i,\Lambda_i$, $i=1,2,3$ and let $R_{ij}(\lambda)=
R_{\Lambda_i,\Lambda_j}(z_i-z_j,\lambda)$.
The dynamical Yang--Baxter equation is equivalent
to the equation
\be
R_{23}^{*(23)}(\lambda-2\eta h^{(1)})
R_{13}^{*(13)}(\lambda)
R_{12}^{*(12)}(\lambda-2\eta h^{(3)})
=
R_{12}^{*(12)}(\lambda)
R_{13}^{*(13)}(\lambda-2\eta h^{(2)})
R_{23}^{*(23)}(\lambda)
\ee
for the dual maps $R^*_{ij}(\lambda)$.
The latter equation follows from using the previous lemma
several times to express $\omega(z_3,z_2,z_1,\lambda)$
in terms of $\omega(z_1,z_2,z_3,\lambda)$ in two different
ways.
The second claim follows by expressing $\omega(w,z,\lambda)$
in terms of $\omega(z,w,\lambda)$ and conversely.
\end{proof}

Let us now consider the case of positive integer weights. In
this case the $R$-matrices have invariant subspaces.
If $\Lambda\in\Z_{\geq 0}$ we let $SV_\Lambda$ be the subspace
of $V_\Lambda$ spanned by $e_{\Lambda+1},\, e_{\Lambda+2},\dots$.
The $\Lambda+1$-dimensional quotient $V_\Lambda/SV_\Lambda$ will
be denoted by $L_\Lambda$, and will be often identified with
$\oplus_{j=0}^m\C e_j$.

\begin{thm}\label{tfd}
Let $z,\eta,\lambda$ be generic and $\Lambda$, $\Mu\in\C$.
\begin{enumerate}
\item[(i)] If $\Lambda\in\Z_{\geq 0}$, then
$R_{\Lambda,\Mu}(z,\lambda)$ preserves $SV_\Lambda\otimes V_\Mu$
\item[(ii)] If $\Mu\in\Z_{\geq 0}$, then
$R_{\Lambda,\Mu}(z,\lambda)$ preserves $V_\Lambda\otimes SV_\Mu$
\item[(iii)] If $\Lambda\in\Z_{\geq 0}$ and $\Mu\in
\Z_{\geq 0}$, then
$R_{\Lambda,\Mu}(z,\lambda)$ preserves $SV_\Lambda\otimes V_\Mu
+V_\Lambda\otimes SV_\Mu$.
\end{enumerate}
\end{thm}

\begin{proof} The claims of this theorem are equivalent
to the statements that the operators dual to the
$R$-matrices preserve the subspaces $L^*_\Lambda\otimes V^*_\Mu$,
$V^*_\Lambda\otimes L^*_\Mu$, $L^*_\Lambda\otimes L^*_\Mu$,
respectively. This follows from Theorem \ref{treso} below,
which gives a characterization of the images by $\omega(z,w,\lambda)$
of these subspaces in terms of residue conditions.
\end{proof}

In particular, if $\Lambda$ and/or $\Mu$ are nonnegative integers,
$R_{\Lambda,\Mu}(z,\lambda)$ induces operators, still denoted
by $R_{\Lambda,\Mu}(z,\lambda)$, on the quotients
$L_\Lambda\otimes V_\Mu$, $V_\Lambda\otimes L_\Mu$ and/or
$L_\Lambda\otimes L_\Mu$. They obey the dynamical Yang--Baxter
equation.

\subsection{Examples}\label{ssexa}
We give now some examples of computations of matrix elements of
the $R$-matrix $R_{\Lambda,\Mu}(z-w,\lambda)$, assuming that the
parameters are generic.

The $R$-matrix is calculated as the transition matrix relating
two bases of $F_{ab}(z,w,\lambda)$: let
\be \tilde\omega_{ij}=
\omega(w,z,\lambda)\,e_j^*\otimes e_i^*, \qquad
\omega_{ij}=\omega(z,w,\lambda)\,e_i^*\otimes e_j^*.
\ee
The matrix elements of $R$ with respect to the basis $e_j\otimes e_k$
are given by $\tilde\omega_{kl}=\sum_{ij}R_{ij}^{kl}\omega_{ij}$.
By construction,
this matrix commutes with $h^{(1)}+h^{(2)}$, and thus preserves
the weight spaces
\be
(V_\Lambda\otimes V_{\Mu})[\Lambda+\Mu-2m]=\oplus_{j=0}^m\C e_j\otimes
e_{m-j}
\ee
We may therefore consider the problem of computing the matrix elements
of the $R$-matrix separately on each weight space.
Without loss of generality we assume that $w=0$.

Let $m=0$. Then the weight space is spanned by $e_0\otimes e_0$ and
$\omega_{00}=\tilde\omega_{00}=1$. Therefore $R_{00}^{00}=1$.

Let now $m=1$. The basis elements are functions of one variable
$t=t_1$ and we have (with $a=\eta\Lambda$, $b=\eta\Mu$)
\be
\omega_{01}(t)=\frac{\theta(\lambda+2\eta+t-2a-b)\theta(t-z+a)}
{\theta(t-b)\theta(t-z-a)},
\qquad
\omega_{10}(t)=\frac{\theta(\lambda+2\eta+t-z-a)}
{\theta(t-z-a)},
\ee
and
\be
\tilde\omega_{01}(t)=\frac{\theta(\lambda+2\eta+t-b)}
{\theta(t-b)}
,
\qquad
\tilde\omega_{10}(t)=\frac{\theta(\lambda+2\eta+t-z-2b-a)\theta(t+b)}
{\theta(t-z-a)\theta(t-b)}
.
\ee
To express one basis in terms of the other we notice that these
functions
have simple poles at two points (modulo the lattice).
The elements of the $R$-matrix are
determined by comparing the residues at these two points.
We have
\bea
\res_{t=b}\,
\omega_{01}(t)=
\frac{\theta(\lambda\!+\!2\eta\!-\!2a)\theta(z\!-\!b\!-\!a)}
{\theta'(0)\theta(z\!-\!b\!+\!a)}
,
& &
\res_{t=z+a}\,\omega_{01}(t)=\frac{\theta(\lambda\!+\!2\eta\!+\!z\!-\!a\!-\!b)
\theta(2a)}
{\theta'(0)\theta(z\!+\!a\!-\!b)}
,\\
\res_{t=b}\,\omega_{10}(t)=0
,
& &
\res_{t=z\!+\!a}\,\omega_{10}(t)=\frac{\theta(\lambda\!+\!2\eta)}
{\theta'(0)}.
\eea
On the other hand,
\bea
\res_{t=b}\,\tilde\omega_{01}(t)=
\frac
{\theta(\lambda\!+\!2\eta)}
{\theta'(0)}
,
& &
\res_{t=z+a}\,\tilde\omega_{01}(t)=0
,
\\
\res_{t=b}\,\tilde\omega_{10}(t)=
\frac
{\theta(\lambda\!+\!2\eta\!-\!z\!-\!b\!-\!a)\theta(2b)}
{\theta'(0)\theta(\!-\!z\!+\!b\!-\!a)}
,
& &
\res_{t=z+a}\,\tilde\omega_{10}(t)=
\frac
{\theta(\lambda\!+\!2\eta\!-\!2b)\theta(z\!+\!a\!+\!b)}
{\theta'(0)\theta(z\!+\!a\!-\!b)}.
\eea
If we put these residues into matrices
$\left(\begin{matrix}A&B\\ 0&D\end{matrix}\right)$ and
$\left(\begin{matrix}\tilde A& 0\\ \tilde C &\tilde
D\end{matrix}\right)$,
respectively, we get that the restriction of $R_{\Lambda,\Mu}(z,\lambda)$
to $(V_\Lambda\otimes V_\Mu)[\Lambda+\Mu-2]$ in the
basis $e_0\otimes e_1$, $e_1\otimes e_0$ is given by the
matrix
\begin{displaymath}
\left(\begin{matrix}
\tilde A&0\\\tilde C&\tilde D
\end{matrix}\right)
\left(\begin{matrix}
A&B\\ 0&D
\end{matrix}\right)^{-1}
=
\left(\begin{matrix}
\tilde AA^{-1}&-\tilde ABA^{-1}D^{-1}\\
\tilde CA^{-1}&-\tilde CBA^{-1}D^{-1}+\tilde DD^{-1}
\end{matrix}\right)
.
\end{displaymath}
These matrix elements have product form except the last one, which
can be further simplified:
\bea
R^{10}_{10}&=&
-\tilde CBA^{-1}D^{-1}+\tilde DD^{-1}\\
&=&
\frac{
\theta(\lambda+2\eta-z-a-b)
\theta(2b)
\theta(\lambda+2\eta+z-a-b)
\theta(2a)}
{\theta(\lambda+2\eta)\theta(z-a-b)
\theta(\lambda+2\eta-2a)
\theta(z+a-b)
}\\
& &+
\frac{\theta(\lambda+2\eta-2b)
\theta(z+a+b)}
{\theta(\lambda+2\eta)\theta(z+a-b)}.
\eea
To simplify this expression, consider it as a function of $\lambda$:
it is a 1-periodic function, which gets multiplied
by $e^{4\pi ib}$ when $\lambda$ is replaced by $\lambda+\tau$.
It is regular at the apparent singularity $\lambda=-2\eta$. So it
has only simple poles at $\lambda=2\eta-2a$ and its translates
by $\Z+\tau\Z$. These properties uniquely define this function
up to a multiplicative function: it can therefore be written in
the form $C\theta(\lambda+2\eta-2a-2b)/\theta(\lambda+2\eta-2a)$.
The constant $C$ can be determined by computing the residue at
$\lambda=2a-2\eta$.

Let us summarize the results of these calculations.

\begin{proposition}\label{Rijkl}
Let the matrix elements of the $R$-matrix $R_{\Lambda,\Mu}(z,\lambda)$
be defined by
\be
R_{\Lambda,\Mu}(z,\lambda)e_i\otimes e_j=\tsize\sum_{kl}
R_{ij}^{kl}
e_k\otimes e_l.
\ee
Then
\bea
R^{00}_{00}&=&1\\
R^{01}_{01}&=&
\frac{\theta(z+\eta\Lambda-\eta\Mu)\theta(\lambda+2\eta)}
{\theta(z-\eta\Lambda-\eta\Mu)\theta(\lambda+2\eta(1-\Lambda))}
\\
R^{01}_{10}
&=&
-
\frac{\theta(\lambda+2\eta+z-\eta\Lambda-\eta\Mu)\theta(2\eta\Lambda)}
{\theta(z-\eta\Lambda-\eta\Mu)\theta(\lambda+2\eta(1-\Lambda))}
\\
R^{10}_{01}
&=&
-\frac{\theta(\lambda+2\eta-z-\eta\Lambda-\eta\Mu)\theta(2\eta\Mu)}
{\theta(z-\eta\Lambda-\eta\Mu)\theta(\lambda+2\eta(1-\Lambda))},
\\
R^{10}_{10}
&=&
\frac{\theta(z+\eta\Mu-\eta\Lambda)\theta(\lambda+2\eta(1-\Lambda-\Mu))}
{\theta(z-\eta\Lambda-\eta\Mu)\theta(\lambda+2\eta(1-\Lambda))}.
\eea
\end{proposition}
For higher values of $m$, the matrix elements can still be computed
along the same lines. The residues are replaced by suitable iterated
residues defined for functions of several variables. This technique
is explained in Section \ref{sp}. The iterated residues of basis
elements $\omega_{ij}$ and $\tilde\omega_{ij}$ form triangular
matrices whose diagonal matrix elements do not vanish for
generic parameters. This gives an expression of the restriction of
the $R$-matrix to any weight space as a product of a lower triangular
matrix by an upper triangular matrix, whose entries are
iterated residues, just as in the one-variable case.

In the next subsection we present the results of this calculation
relevant for the representation theory of elliptic quantum groups.

\subsection{Evaluation Verma modules and their tensor products}
\label{ssevm}
Here we explain the relation between the geometric construction of tensor
products and $R$-matrices and the representation theory
of $E_{\tau,\eta}(sl_2)$ \cite{FV1}.

We first recall the definition of a representation of $E_{\tau,\eta}(sl_2)$:
let $\h$ act on $\C^2$ via $h={\mathrm {diag}}(1,-1)$. A representation
of $E_{\tau,\eta}(sl_2)$ is
an $\h$-module $W$ with diagonalizable action of $h$ and finite
dimensional eigenspaces, together with an operator $L(z,\lambda)
\in\End(\C^2\otimes W)$ (the ``$L$-operator''),
commuting with $h^{(1)}+h^{(2)}$, and
obeying the relations
\bea
R^{(12)}(z\!-\!w,\lambda\!-\!2\eta h^{(3)})\!\!\!\!\!
&L^{(13)}(z,\lambda)\,
L^{(23)}(w,\lambda-2\eta h^{(1)})&\!\!\!\!\!
\\
&=\,L^{(23)}(w,\lambda)\,
L^{(13)}(z,\lambda\!-\!2\eta h^{(2)})&\!\!\!\!\!
R^{(12)}(z\!-\!w,\lambda)
\eea
in $\End(\C^2\otimes\C^2\otimes W)$. The
{\em fundamental $R$-matrix} $R(z,\lambda)\in\End(\C^2\otimes\C^2)$
is the following solution of the
dynamical Yang--Baxter equation:
let $e_0$, $e_1$ be the standard basis of $\C^2$, then
with respect to the basis $e_0\otimes e_0$, $e_0\otimes e_1$,
$e_1\otimes e_0$, $e_1\otimes e_1$ of $\C^2\otimes\C^2$,
\be
R(z,\lambda)=
\left(\begin{matrix}
1 & 0 & 0 & 0\\
0 & \alpha(z,\lambda) & \beta(z,\lambda) & 0\\
0 & \beta(z,-\lambda) & \alpha(z,-\lambda) & 0\\
0 & 0 & 0 & 1
\end{matrix}\right).
\ee
where
\be
\alpha(z,\lambda)=
\frac{\theta(\lambda+2\eta)\theta(z)}
{\theta(\lambda)\theta(z-2\eta)},\qquad
\beta(z,\lambda)=
-\frac{\theta(\lambda+z)\theta(2\eta)}
{\theta(\lambda)\theta(z-2\eta)}.
\ee
To discuss representation theory, it is convenient
to think of $L(z,\lambda)\in\End(\C^2\otimes V)$ as a two by two
matrix with entries $a(z,\lambda)$, $b(z,\lambda)$, $c(z,\lambda)$,
$d(z,\lambda)$ in $\End(W)$. In \cite{FV1} we wrote explicitly the
relations that these four operators must satisfy, and defined
a class of representations, the evaluation Verma modules,
by giving explicitly the action of these four operators on
basis vectors. Next, we show how these formulae can be
obtained from our geometric construction.

\begin{thm}\label{tfrm} Let us identify $L_1$ with $\C^2$ via the basis
$e_0,e_1$.
The $R$-matrix $R_{1,1}(z,\lambda)\in\End(L_1\otimes L_1)$
coincides with the fundamental $R$-matrix.
\end{thm}
\begin{proof} Most of the matrix elements have been computed in
\ref{ssexa}. The only remaining element to compute is $R^{11}_{11}$.
By Theorem \ref{tfd}, we know that the dual of the $R$-matrix
preserves the one-dimensional space spanned by $e_1^*\otimes e_1^*$
if $\Lambda=\Mu=1$.
This means that $\omega_{1,1}=\omega(z,w,\lambda)e_1^*\otimes e_1^*$
is proportional to $\tilde\omega_{1,1}=\omega(w,z,\lambda)e_1^*\otimes e_1^*$.
We have, by definition,
\bea
\omega_{1,1}&=&
\frac{\theta(\lambda+t_1+\eta-z)\theta(\lambda+t_2+3\eta-w)\theta(t_2-z+\eta)}
{\theta(t_1-z-\eta)\theta(t_2-w-\eta)\theta(t_2-z-\eta)}\\
& &+
\frac{
\theta(\lambda+t_2+\eta -z)
\theta(\lambda+t_1+3\eta-w)
\theta(t_1-z+\eta)\theta(t_1-t_2-2\eta)}
{\theta(t_2-z-\eta)\theta(t_1-w-\eta)\theta(t_1-z-\eta)\theta(t_1-t_2+2\eta)}.
\eea
The formula for $\tilde\omega_{1,1}$ is obtained by interchanging $z$
and $w$.
To find the constant of proportionality (the inverse of $R_{11}^{11}$),
we compute both expressions at some point. We can take
for instance $t_1=2\eta$, $t_2=0$. Then it is obvious that at this
point $\omega_{1,1}
=\tilde\omega_{1,1}$. Thus $R_{11}^{11}=1$.
\end{proof}

\begin{corollary}
For any $w,\Mu\in\C$, the $\h$-module $V_\Mu$ together with the
operator $L(z,\lambda)=R_{1,\Mu}(z-w,\lambda)\in\End(L_1\otimes V_\Mu)$
defines a representation of $E_{\tau,\eta}(sl_2)$.
\end{corollary}

This representation is called in \cite{FV1}
the evaluation Verma module with evaluation
point $w$ and highest weight $\Mu$. It
is denoted by $V_\Mu(w)$.
The matrix elements of $L(z,\lambda)$ can be computed
straightforwardly by the method described in the previous section.
The result is given explicitly in \cite{FV1}, Theorem 3, in terms
of the action of $a(z,\lambda)$,\dots, $d(z,\lambda)$. In the notation
of Proposition \ref{Rijkl}, this result amounts to the following
formulae for the matrix elements $R_{ij}^{kl}$ of $R_{1,\Lambda}(z,\lambda)\in
\End(L_1\otimes V_\Lambda)$.
\bea
R_{0k}^{0k}&=&{\theta(z-(\Lambda+1-2k)\eta)\over
\theta(z-(\Lambda+1)\eta)} {
\theta(\lambda+2k\eta) \over
\theta(\lambda) },
\\
R_{1k}^{0,k+1}&=&-\,
{\theta(\lambda+z-(\Lambda-1-2k)\eta)\over
\theta(z-(\Lambda+1)\eta)}{
\theta(2\eta) \over \theta(\lambda) }, \\
R_{0k}^{1,k-1}&=&-\,{\theta(\lambda-z-(\Lambda+1-2k)\eta)\over
\theta(z-(\Lambda+1)\eta)} {
\theta(2(\Lambda+1-k)\eta)\over
\theta(\lambda) } {\theta(2k\eta)\over \theta(2\eta)}, \\
R_{1k}^{1k}&=&
{\theta(z-(-\Lambda+1+2k)\eta)\over
\theta(z-(\Lambda+1)\eta) }{
\theta(\lambda-2(\Lambda-k)\eta) \over
\theta(\lambda) }.
\eea
These formulae can be obtained by computing residues, as
explained in the previous subsection.
Moreover, the tensor product construction of \ref{sstp} is
related to the tensor product of representations of the elliptic
quantum group. Recall that if $W_1$, $W_2$ are representations
of the elliptic quantum group
with $L$-operators $L_1(z,\lambda)$, $L_2(z,\lambda)$, then
their tensor product $W=W_1\otimes W_2$ with $L$-operator
\be
L(z,\lambda)=L_1(z,\lambda-2\eta h^{(3)})^{(12)}
L_2(z,\lambda)^{(13)}\in\End(\C^2\otimes W)
\ee
is also a representation of the elliptic quantum group.

\begin{thm} Let $\Lambda_1,\dots,\Lambda_n\in\C$ and
$z_1,\dots,z_n$ be generic complex numbers.
Let $V=V_{\Lambda_1}\otimes\cdots\otimes V_{\Lambda_n}$ and
$L(z,\lambda)\in\End(V_{\Lambda=1}\otimes V)$ be defined
by the relation
\be
\omega(z,z_1,\dots,z_n,\lambda)L(z,\lambda)^*=
\omega(z_1,\dots,z_n,z,\lambda)P
\ee
in $\End((V_1\otimes V)^*)=\End(V_1^*\otimes V^*)$,
where $Pv_1\otimes v=v\otimes v_1$, if $v_1\in V_1^*$, $v\in V^*$.
Then $L(z,\lambda)$ is well-defined as an endomorphism of the quotient
$L_1\otimes V=\C^2\otimes V$, and defines the structure of
a representation of $E_{\tau,\eta}(sl_2)$ on $V$.
This representation is isomorphic to the tensor product
of evaluation Verma modules
\be
V_{\Lambda_n}(z_n)\otimes\cdots\otimes V_{\Lambda_1}(z_1),
\ee
with isomorphism $u_1\otimes\cdots\otimes u_n\mapsto
u_n\otimes\cdots\otimes u_1$
\end{thm}

\begin{proof}
We prove this if $n=2$. The proof in the general case is
the same, but requires more writing.

By applying Lemma \ref{lsym} twice, we get
\begin{gather*}
\omega(z,z_1,z_2,\lambda)R_{1,\Lambda_1}^*(z-z_1,\lambda)^{(12)}
R_{1,\Lambda_2}^*(z-z_2,\lambda-2\eta h^{(2)})^{(13)}
\\
=
\omega(z_1,z,z_2,\lambda)
R_{1,\Lambda_2}^*(z-z_2,\lambda-2\eta h^{(1)})^{(23)}P^{(12)}
\\
=\omega(z_1,z_2,z,\lambda)P^{(23)}P^{(12)}.
\end{gather*}
But $P^{(23)}P^{(12)}=P$, the flip $V_1\otimes V\to V\otimes V_1$.
By taking duals, we obtain
\be
L(z,\lambda)=R_{1,\Lambda_2}(z-z_2,\lambda-2\eta h^{(2)})^{(13)}
R_{1,\Lambda_1}(z-z_1,\lambda)^{(12)}.
\ee
This operator in well-defined as an endomorphism of
the quotient $L_1\otimes V$, since its factors are. Moreover
\be
P^{(23)}L(z,\lambda)P^{(23)}
=R_{1,\Lambda_2}(z-z_2,\lambda-2\eta h^{(3)})^{(12)}
R_{1,\Lambda_1}(z-z_1,\lambda)^{(13)},
\ee
which is the $L$-operator of the representation $V_{\Lambda_2}(z_2)
\otimes V_{\Lambda_1}(z_1)$.
\end{proof}

Finally, the dynamical Yang--Baxter equation in $L_1\otimes V_\Lambda
\otimes V_\Mu$ can be stated as saying that
$R_{\Lambda,\Mu}(z-w,\lambda)P$ is
an isomorphism from $V_\Mu(w)\otimes V_\Lambda(z)$ to
$V_\Lambda(z)\otimes V_\Mu(w)$, see \cite{FV1}. Therefore, by uniqueness,
the $R$-matrices constructed here coincide with the solutions
of the dynamical Yang--Baxter equation described in Section 13
of \cite{FV1}.

\section{Bethe ansatz}\label{sba}

In this section we fix the parameters $\eta$, $\tau$ of the
elliptic quantum group and $n$ complex numbers
$\Lambda_1,\dots,\Lambda_n$ such that $\sum \Lambda_l
=2 m$ for some nonnegative integer $m$.
We also set $a_j=\eta\Lambda_j$, $j=1,\dots,n$.
Let $z_1,\dots,z_n$ be generic complex numbers.
Then the zero
weight space $V[0]$ of $V=V_{\Lambda_1}\otimes\cdots\otimes V_{\Lambda_n}$
is non-trivial, and we have commuting difference operators
$H_j=H_j(z_1,\dots,z_n)$, $j=1,\dots,n$, cf.\ \Ref{Go},
acting on the space Fun$(V[0])$.

We give a formula for common quasiperiodic eigenfunctions of the operators
$H_j$, i.e., functions $\psi$ such that $H_j\psi=\epsilon_j\psi$,
$j=1,\dots, n$, for some $\epsilon_j\in\C$. The quasiperiodicity
assumption means that we require that $\psi(\lambda+1)=
\mu\psi(\lambda)$ for some multiplier $\mu\in\C^\times$.

Our formula is given in terms of the expressions of the
previous section. For $j_1,\dots,j_n\in\Z_{\geq 0}$ such that
$\sum j_i=m$, let
$\omega_{j_1,\dots, j_n}(t_1,\dots,t_m,\lambda)$ be
the value at $t_1,\dots,t_m$ of the function
$\omega(z_1,\dots,z_n,\lambda)
e_{j_1}^*\otimes\cdots\otimes e_{j_n}^*$.
\begin{thm}\label{tba}
(cf.\ \cite{TV1})
Let $c\in\C$.
Suppose that $t_1,\dots,t_m$ obey the system of
``Bethe ansatz'' equations
\begin{equation}\label{bae}
\prod_{l=1}^{n}
\frac{\theta(t_j-z_l+a_l)}
{\theta(t_j-z_l-a_l)}
\prod_{k:k\neq j}
\frac{\theta(t_j-t_k-2\eta)}{\theta(t_j-t_k+2\eta)}
=e^{-4\eta c},\qquad j=1,\dots,m.
\end{equation}
Then
\begin{equation}\label{psi}
\psi(\lambda)=\sum_{j_1+\cdots+j_n=m}
e^{c\lambda}
\omega_{j_1,\dots,j_n}(t_1,\dots,t_m,\lambda)\,
e_{j_1}\otimes\cdots\otimes e_{j_n}
\end{equation}
is a common eigenfunction of the commuting operators
$H_j$, $j=1,\dots,n$, with eigenvalues
\be
\epsilon_j=e^{-2ca_j}\prod_{k=1}^m
\frac{\theta(t_k-z_j-a_j)}
{\theta(t_k-z_j+a_j)}
\ee
and multiplier $\mu=(-1)^me^c$. Moreover, if $\,t_1,\dots,t_n$ are a
solution of
\Ref{bae}, and $\sigma\in S_m$ is any permutation,
then $t_{\sigma(1)},\dots,t_{\sigma(n)}$ are also a solution.
The eigenfunctions $\psi$ corresponding to these two solutions
are proportional.
\end{thm}

\begin{example}
Let $m=1$. Thus $\Lambda_1+\cdots+\Lambda_n=2$.
The Bethe ansatz equation for $t=t_1$ and $c$ is
\be
\prod_{l=1}^n
\frac{\theta(t-z_l-\eta\Lambda_l)}
{\theta(t-z_l+\eta\Lambda_l)}
=e^{4\eta c}.
\ee
If $t,c$ obey this equation, then
\be
\psi(\lambda)=e^{c\lambda}
\sum_{l=1}^n
\left(
\frac
{\theta(\lambda\!+\!t\!-\!z_l\!-
\!\eta\Lambda_l\!-\!2\eta\sum_{j=1}^{l-1}{\Lambda_j})}
{\theta(t\!-\!z_l\!-\!\eta\Lambda_l)}
\prod_{j=1}^{l-1}
\frac
{\theta(t\!-\!z_j\!+\!\eta\Lambda_l)}
{\theta(t\!-\!z_j\!-\!\eta\Lambda_l)}\right)
e_0\otimes\cdots\otimes e_1\otimes\cdots\otimes e_0,
\ee
($e_1$ is at the $l$th position) is a common eigenfunction of the
difference operators $H_j(z)$ $j=1,\dots,n$. The eigenvalues are
\be
\epsilon_j=e^{-2\eta c\Lambda_j}\frac{\theta(t-z_j-\eta\Lambda_j)}
{\theta(t-z_j+\eta\Lambda_j)}.
\ee
\end{example}

The proof of this theorem occupies the rest of this section.
The fact that solutions are mapped to solutions by permutations
is clear from the form of the equations \Ref{bae}.
The functions $\omega_{j_1,\dots,j_n}$ are symmetric under
the action of permutations of Lemma \ref{lsm}. Therefore
functions $\psi(\lambda)$ corresponding to solutions
related by permutations are proportional.

We are left to prove that $\psi$ is an eigenfunction of the
operators $H_j$.

We first reduce the problem to showing that
$\psi(\lambda)$ is an eigenfunction of $H_1$ for all permutations
of the parameters $z_i,\Lambda_i$.
For this we must study the dependence
of $H_j$ on the ordering of the parameters $z_i,\Lambda_i$.
Let us introduce some notation to indicate explicitly
the dependence on the parameters: let $S_m$ act on $\C^m$ by
permutations of the coordinates, and let $s_j$ be the
transposition of the $j$th and $(j+1)$st coordinates. Let
us write $H_j({ \Lambda},{ z})$ for the operator $H_j$
acting on the space of functions $\cal H({ \Lambda})$ of
functions of $\lambda$ with values in $(V_{\Lambda_1}
\otimes\cdots\otimes V_{\Lambda_n})[0]$,
${ \Lambda},{ z}\in\C^n$.

\begin{lemma}\label{lord} Let ${ \Lambda}=(\Lambda_1,
\dots,\Lambda_n)$, ${ z}=
(z_1,\dots,z_n)\in\C^n$ and let $S_j:{\cal H}
({ \Lambda})\to {\cal H}(s_j{ \Lambda})$
$(j=1,\dots n-1)$ be the linear map
\be
S_jf(\lambda)=P^{(j,j+1)}R_{\Lambda_j\Lambda_{j+1}}
^{(j,j+1)}(z_j\!-\! z_{j+1},\lambda
-2\eta(h^{(1)}+\cdots+h^{(j-1)}))f(\lambda).
\ee
Then $H_{j+1}(s_j{ \Lambda},
s_j{ z})S_j=S_jH_j({ \Lambda},{ z})$,
for all $j=1,\dots, n-1$.
\end{lemma}
\begin{proof}
The proof is a matter of writing the claim using the
representation of $H_j$ given in \Ref{Go}, and
bringing $P^{(j,j+1)}$ to the left of the $R$-matrices,
using commutation relations.
\end{proof}

Notice that the equations \Ref{bae} do
not depend on the ordering of the parameters $z_i,\Lambda_i$. Let
us fix a solution $t^*$ of \Ref{bae} with parameters
${\Lambda}=(\Lambda_1,\dots,\Lambda_n)$
and ${z}=(z_1,\dots,z_n)$. Then
for each permutation $\sigma\in S_n$, $t^*$ is still a
solution of the equations \Ref{bae} with parameters $\sigma{\Lambda}$,
$\sigma{z}$ and we have a corresponding function
$\psi_\sigma(\lambda)\in {\cal H}(\sigma{\Lambda})$
given by the formula \Ref{psi} with
parameters $\sigma{\Lambda}$, $\sigma{z}$. It takes
values in $V_{\Lambda_{\sigma^{-1}(1)}}
\otimes\cdots\otimes V_{\Lambda_{\sigma^{-1}(n)}}$.
We may now rephrase Lemma \ref{lsym} as follows:
\be
\psi_{s_j}(\lambda)=P^{(j,j+1)}
R_{\Lambda_j\Lambda_{j+1}}
^{(j,j+1)}(z_j-z_{j+1},\lambda-2\eta(h^{(1)}+\cdots+h^{(j-1)}))\psi(\lambda).
\ee
With Lemma \ref{lord}, we see that the fact that $\psi(\lambda)$
is an eigenfunction of $H_j({\Lambda},{z})$ with eigenvalue
$\epsilon_j=\epsilon_j({\Lambda},{z})$
for all $j=1,\dots,n$ is equivalent to the
fact that $\psi_\sigma(\lambda)$ is an eigenfunction
of $H_1(\sigma{\Lambda},\sigma{z})$ with eigenvalue
$\epsilon_1(\sigma{\Lambda},\sigma{z})$ for all $\sigma\in S_n$.

Therefore, it is sufficient to prove that $\psi$ is an eigenfunction
of $H_1$ with eigenvalue $\epsilon_1$ for any solution of
\Ref{bae} with arbitrary parameters.

\begin{lemma}\label{l13}
Let $u\in V_{\Lambda_1}^*[\mu]$, $v\in (V_{\Lambda_2}^* \otimes\cdots\otimes
V_{\Lambda_{n}}^*)[-\mu]$. Thus $\mu=\Lambda_1-2m'=
-\sum_{j=2}^{n}\Lambda_j+2m''$ for some integers $m', m''\geq 0$. Set
$m=m'+m''$. Let $f=\omega(z_1,\lambda)u$,
$g=\omega(z_2,\dots,z_{n},\lambda-2\eta\mu)v$. Then
$\omega(z_1,\dots,z_n,\lambda)$ maps $u\otimes v$ to
\be
\frac1{(n-1)!}\Sym\biggl( f(t_1,\dots,t_{m'}) g(t_{m'+1},\dots,t_{m})
\prod_{j>m'} \frac{\theta(t_j-z_1+a_1)} {\theta(t_j-z_1-a_1)}\biggr),
\ee
and $\omega(z_2,\dots,z_{n},z_1,\lambda-2\eta\mu)$ maps $v\otimes u$
to
\be
\frac1{(n\!-\!1)!}\Sym\biggl( g(t_{m'+1},\dots,t_{m})
f(t_1,\dots,t_{m'}) \prod_{\begin{matrix}{\scriptstyle{j\leq m'}}\\
\scriptstyle{l\geq2}\end{matrix}} \frac{\theta(t_j\!-\!z_l\!+\!a_l)}
{\theta(t_j\!-\!z_l\!-\!a_l)} \prod_{\begin{matrix}{\scriptstyle{j\leq
m'}}\\ \scriptstyle{k>m'}\end{matrix}} \frac
{\theta(t_j\!-\!t_k\!-\!2\eta)} {\theta(t_j\!-\!t_k\!+\!2\eta)} \biggr).
\ee
\end{lemma}

\begin{proof}
The first formula is a rewriting of the definition. In the
second formula, the expression appearing in the definition was
replaced by another term in the sum over $S_m$ defining $\Sym$,
which gives the same result.
\end{proof}

\begin{lemma}\label{l14} If $t^*=(t^*_1,\dots,t^*_m)$ is a solution
of \Ref{bae} and $u$, $v$ are as in
Lemma \ref{l13}, then
\be
\omega(z_2,\dots,z_{n},z_1,\lambda-2\eta\mu)
v\otimes u|_{t^*}=
e^{2c\eta\mu}
\epsilon_1\cdot\omega(z_1,\dots,z_n,\lambda)u\otimes v|_{t^*}
\ee
Here, ${}|_{t^*}$ denotes the value at $t^*$ of a function
of $t_1,\dots,t_m$.
\end{lemma}
\begin{proof}
Let us introduce the ratio $h(t_1,\dots,t_m)$ of the
factors appearing in the previous lemma:
\be\prod_{\begin{matrix}{\scriptstyle{j\leq m'}}\\
\scriptstyle{l>2}\end{matrix}}
\frac{\theta(t_j\!-\!z_l\!+\!a_l)}
{\theta(t_j\!-\!z_l\!-\!a_l)}
\prod_{\begin{matrix}{\scriptstyle{j\leq m'}}\\
\scriptstyle{k>m'}\end{matrix}}
\frac
{\theta(t_j\!-\!t_k\!-\!2\eta)}
{\theta(t_j\!-\!t_k\!+\!2\eta)}
=
h(t_1,\dots,t_m)
\prod_{j>m'}
\frac{\theta(t_j-z_1+a_1)}
{\theta(t_j-z_1-a_1)}\, .
\ee
This equation can be rewritten as
\be
h(t_1,\dots,t_m)
=
e^{2ca_1}
\epsilon_1(t_1,\dots,t_m)
\prod_{j\leq m'}\biggl(
\prod_{l=1}^n
\frac{\theta(t_j\!-\!z_l\!-\!a_l)}
{\theta(t_j\!-\!z_l\!+\!a_l)}
\prod_{k>m'}
\frac
{\theta(t_j\!-\!t_k\!+\!2\eta)}
{\theta(t_j\!-\!t_k\!-\!2\eta)}\biggr),
\ee
where $\epsilon_1(t_1\dots,t_m)=e^{-2ca_1}
\prod_{j=1}^m
{\theta(t_j-z_1-a_1)}/{\theta(t_j-z_1+a_1)}$.
The product of the first $m'$ Bethe ansatz equations
\Ref{bae}, implies that, for $t=t^*$,
\be
h(t_1,\dots,t_m)
=e^{2ca_1}e^{-4\eta cm'}\epsilon_1(t_1,\dots,t_m)
\ee
The right-hand side of this equation is symmetric under
permutations of $t_1,\dots,t_m$. Moreover,
if $t_1,\dots,t_m$ is a solution of \Ref{bae}
then also $t_{\sigma(1)},\dots, t_{\sigma(m)}$ for any permutation
$\sigma$. We conclude that for all solutions of \Ref{bae},
\be
h(t_1,\dots,t_m)=h(t_{\sigma(1)},\dots,t_{\sigma(m)}),
\ee
for all $\sigma\in S_m$. It follows that terms corresponding
to the same permutation in the two equations of Lemma
\ref{l13} are proportional, with the same constant of
proportionality $e^{2c\eta\mu}\epsilon_1$.
\end{proof}

Now we can prove Theorem \ref{tba}: $H_1$ has
the form $H_1f(\lambda)=\Gamma_1\tilde H_1f(\lambda)$, where
$\tilde H_1$ is the operator of multiplication by
\begin{equation}\label{tH}
\tilde H_1(\lambda)=
R_{\Lambda_1,\Lambda_n}(z_1\!-\!z_n,\lambda-2\eta\tsize\Sum_{l=2}^{n-1}h^{(l)})
^{(1n)}\cdots
R_{\Lambda_1,\Lambda_2}(z_1\!-\!z_2,\lambda)^{(12)}.
\end{equation}
By using $n\!-\!1$ times Lemma \ref{lsym}, we see that
\be
\omega(z_2,\dots,z_n,z_1,\lambda)
P^{(n,n-1)}\cdots P^{(3,2)}P^{(2,1)}
=\omega(z_1,\dots,z_n,\lambda)
\tilde H_1(\lambda)^*.
\ee
Now we use Lemma \ref{l14}. Let $t^*$ be a solution
of \Ref{bae}. Then, for any $u\in V_{\Lambda_1}^*[\mu]$, $v\in (V_{\Lambda_2}^*
\otimes\cdots\otimes V_{\Lambda_{n}}^*)[-\mu]$,
\begin{equation}\label{e09}
\omega(z_1,\dots,z_n,\lambda-2\eta\mu)\tilde
H_1(\lambda-2\eta\mu)^*u\otimes v|_{t^*}
=e^{2\eta c \mu}\epsilon_1\cdot\omega(z_1,\dots,z_n,\lambda)u\otimes v|_{t^*}.
\end{equation}
The fact that $\psi(\lambda)$ is an eigenfunction is
a consequence of this formula:
the function $\psi\in \mbox{Fun}V[0]$
is uniquely characterized by the property that, for any
$w\in (\otimes_j V_{\Lambda_j}^*)[0]$, the pairing with $\psi(\lambda)$
is
\be
\langle w,\psi(\lambda)\rangle=
e^{c\lambda}\omega(z_1,\dots,z_n,\lambda)w|_{t^*}.
\ee
Take $w=u\otimes v$, as above.
Then $\langle w,H_1\psi(\lambda)\rangle
=\langle w,\tilde H_1(\lambda-2\eta\mu)\psi(\lambda
-2\eta\mu)\rangle
=e^{c\lambda-2\eta c\mu}\omega(z_1,\dots,z_n,\lambda-2\eta\mu)
\tilde H_1(\lambda-2\eta\mu)^*w|_{t^*}$. By
\Ref{e09}, this is equal
to $\epsilon_1\langle w,\psi(\lambda)\rangle$. Since
an arbitrary $w$ can be written as linear combination
of vectors of the form $u\otimes v$ as above,
the theorem is proved.

\section{Completeness of Bethe vectors}\label{scbv}

We assume that $2\eta=1/N$, for some positive odd integer $N$.
We also suppose that $\Lambda_1,\dots,\Lambda_n\in\Z$.
This implies that the operators $h^{(i)}$ on $V=V_{\Lambda_1}
\otimes\cdots\otimes V_{\Lambda_n}$ have integer eigenvalues.
Let $z\in\C^n$ be generic.
The commuting
difference operators $H_j(z)$ have coefficients
which are 1-periodic functions of $\lambda$.
Thus, for each $\alpha\in\C^\times$ we may consider the eigenvalue problem
with multiplier condition
\bean\label{bloch}
H_j(z)\psi&\!\!=\!\!&\epsilon_j\psi, \qquad j=1,\dots,n, \\
\psi(\lambda+1)&\!\!=\!\!&\alpha\psi(\lambda).\notag
\eean
Let $K_N$ be the field of $1/N$-periodic meromorphic functions
of $\lambda\in\C$. The operators $H_j(z)$ are $K_N$-linear
on the $N\dim(V[0])$-dimensional $K_N$-vector space
${\cal H}_\alpha=\{\psi\in\mbox{Fun}(V[0]): \psi(\lambda\!+\!1)=
\alpha\psi(\lambda)\}$.

\begin{thm}\label{toc}
Suppose that $\Lambda_1,\dots,\Lambda_n$ are integers
and let $\sum_{j=1}^n\Lambda=2m$, $m\in\Z_{\geq 0}$.
Let $2\eta N=1$, for some odd integer $N> m$.
For generic $\alpha$, there are $N\dim(V[0])$
solutions of the Bethe ansatz equations \Ref{bae} with
$e^c=(-1)^m\alpha$, such that the corresponding eigenfunctions
form a basis of the $K_N$-vector space ${\cal H}_\alpha$
\end{thm}
\begin{proof}
Let us first consider the case $m=1$, and let $t=t_1$
(see the example after Theorem \ref{tba}).
The Bethe ansatz equation
is
\be
\prod_{l=1}^n
\frac{\theta(t-z_l-\eta\Lambda_l)}
{\theta(t-z_l+\eta\Lambda_l)}
=e^{2c/N},
\ee
corresponding to the eigenfunction
\be
\psi(\lambda)=e^{c\lambda}\sum_{j=1}^n\omega_{(j)}(t,\lambda)
\,e_0\otimes\cdots
\otimes e_1\otimes\cdots\otimes e_0,
\ee
where $\omega_{(j)}$ stands for $\omega_{0,\dots,0,1,0,\dots,0}$
with the 1 at the $j$th position.

We want to fix the multiplier $\alpha$. Thus $c$ is of the form
$c_0+2\pi ir$, $r\in\Z$, where $c_0=\ln(\alpha)+m\pi i$ for some
choice of the branch of the logarithm. Let us now take $\alpha$ real
and tending to infinity. Then, for each $k$, in the vicinity of
the point $t_\infty=z_k-\eta\Lambda_k$, we may introduce the local
coordinate $u=t-t_\infty$, and the Bethe ansatz equation
has the form const $u^{-1}(1+O(u))=\alpha^{2/N}e^{-4\pi is/N}$, for
some integer $s$ related to $r$. This
equation has a solution for any sufficiently large $\alpha$. It has the
asymptotic form $t\sim z-a_k+\mbox{const}\,\alpha^{-2/N}e^{4\pi is/N}$.
Let us denote this solution by
$t^{k,s}=t^{k,s}(\alpha)$. The solutions corresponding to different $k,s$
are distinct if the
indices run over the sets $k=1,\dots,n$, $s=0,\dots, N-1$. We claim that
the corresponding eigenfunctions $\psi^{k,s}$
are linearly independent over
$K_{N}$. We use the following simple observation.

\begin{lemma} If $f_1(\lambda),\dots,f_m(\lambda)$
are meromorphic functions on $\C$ with
values in $\C^d$ that
are linearly dependent over $K_N$, then the vectors
$(f_i(\lambda),f_i(\lambda+1/N),\dots,f_i(\lambda+(N-1)/N))$ in
$\C^{md}$ are linearly dependent over $\C$ for all generic $\lambda\in\C$.
\end{lemma}

\begin{proof}
{}From $\sum_ia_i(\lambda)f_i(\lambda)=0$, with $a_i$ $1/N$-periodic,
we deduce $\sum_ia_i(\lambda)f_i(\lambda+j/N)$, $j=0,\dots,N-1$.
\end{proof}

Therefore it is sufficient to show that the vectors
\be
\tilde\psi^{k,s}(\lambda)=(\psi^{k,s}(\lambda),
\psi^{k,s}(\lambda+1/N),
\dots,
\psi^{k,s}(\lambda+(N-1)/N)) \in V[0]^N\ee
are linearly independent over $\C$ for
generic $\lambda$.

We prove this by showing that the determinant of these vectors does
not vanish for all large enough $\alpha$. The determinant is proportional to
the determinant of the
matrix
\be
A_{(j,r),(k,s)}=
e^{4\pi is(\lambda+r/N)}\omega_{(j)}(t^{k,s},\lambda+r/N),
\qquad j,k=1,\dots,n,\qquad r,s=0,\dots, N-1
\ee
Let now $\alpha$ tend to infinity. Then the solution $t^{k,s}$ approaches
$z_k-a_k$. Thus, $A_{(j,r),(k,s)}$ tends to zero if $k<j$, since
$\omega_j(t=z_k-a_k,\lambda)=0$ for all $\lambda$, owing to the
vanishing of the factor $\theta(t_j-z_k+a_k)$.
Therefore the limiting determinant is block-triangular, and equals
the product of determinants of the diagonal $N$ by $N$ blocks. We have
\be
\det(A)\to
\prod_{j=1}^n\det{}_{r,s}(e^{4\pi is(\lambda+r/N)}
\omega_{(j)}(t=z_j-a_j,\lambda+r/N)), \qquad \alpha\to\infty
\ee
Each factor in the product is, up to multiplication of rows and
columns by nonzero factors, the Vandermonde determinant
\be
\det{}_{0\leq r,s<N}(e^{4\pi irs/N})=\prod_{N>r>s\geq 0}(e^{4\pi ir/N}-e^{4\pi
is/N})
\ee
which does not vanish if $N$ is odd.

The general case is treated along similar lines. We have
solutions $t^{K,s}(\alpha)$ labeled by $K=(k_1,\dots,k_n)$ with $k_j\geq 0$
and $\sum k_i=m$,
and $s\in\{0,\dots,N-1\}$. As the multiplier $\alpha$ tends to infinity,
$t^{K,s}$ converges to the point with coordinates
\be
t_j=z_l-a_l-2\eta(k_1+\cdots+k_l-j),\qquad \mathrm{if}\;
k_1+\cdots+k_{l-1}<j\leq k_1+\cdots+k_l.
\ee
The hypothesis $N> m$ ensures that these coordinates are distinct
modulo the lattice.
By the same arguments as above, one sees that the linear independence
of the corresponding eigenfunctions follows from the non-vanishing
of a determinant. This determinant has for $\alpha\to\infty$
a block-triangular form, and the calculation is reduced to
the calculation of a Vandermonde determinant. The details
are left to the reader.
\end{proof}

Rather than considering the difference operators as acting on
the functions of the continuous parameter $\lambda$, we may
consider the following discrete variant. Let $\mu$ be generic
and $C_\mu=\{\mu+j/N| j\in\Z\}$. Then the commuting difference
operators are defined on the space $\Fun_\mu(V[0])$ of
functions of $\lambda\in C_\mu$ with values in $V[0]$, and preserve
the subspace ${\cal H}_{\alpha,\mu}$ of functions $f$ in $\Fun_\mu(V[0])$
such that $f(\lambda+1)=\alpha f(\lambda)$.

\begin{corollary}\label{cdm}
Suppose that $\Lambda_1,\dots,\Lambda_n$ are integers
and let $\sum_{j=1}^n\Lambda=2m$, $m\in\Z_{\geq 0}$.
Let $2\eta N=1$, for some odd integer $N> m$.
For generic $\alpha$, there are $N\dim(V[0])$
solutions of the Bethe ansatz equations \Ref{bae} with
$e^c=(-1)^m\alpha$, such that the corresponding eigenfunctions
form a basis of the complex vector space ${\cal H}_{\alpha,\mu}$
\end{corollary}

It should be possible to treat the case of general rational $\eta$ (or
more generally in $\Q+\tau\Q$) in a similar way.

\section{Comparison with the algebraic Bethe ansatz}\label{caba}

In this section we compare the Bethe ansatz of this article with
the results of \cite{FV2}, where eigenfunctions of the transfer
matrix of highest weight representations of $E_{\tau,\eta}(sl_2)$
are given.

Let $W$ be a representation of $E_{\tau,\eta}(sl_2)$, see \ref{ssevm}.
Then we have
four operators, $a(z,\lambda)$, $b(z,\lambda)$, $c(z,\lambda)$,
$d(z,\lambda)$, the matrix elements of the $L$-operators, acting
on $W$, and obeying the various relations of $E_{\tau,\eta}(sl_2)$.
The transfer matrix $T_W(z)\in\End(\Fun(W[0]))$ acts on functions
by
\be
T_W(z)f(\lambda)=a(z,\lambda)f(\lambda-2\eta)
+d(z,\lambda)f(\lambda+2\eta).
\ee
The relations imply that $T_W(z)T_W(w)=T_W(w)T_W(z)$ for all $z,w\in\C$.
If $W=L_1(z_n)\otimes\cdots\otimes L_1(z_1)$,
this transfer matrix coincides with the
transfer matrix of \ref{ssrem} conjugated by
\be
\Pi:v_1\otimes\cdots\otimes v_n\to v_n\otimes\cdots\otimes v_1.
\ee
Therefore, in this special case, the
commuting operators $\Pi H_j(z_1,\dots,z_n)\Pi^{-1}$
are equal to $T_W(z_j)$.

In any case it can be shown in general, using the intertwining property of the
$R$-matrices that $T_W(z)$ {\em commutes} with the operators
$\Pi H_j(z_1,\dots,z_n)\Pi^{-1}$
if $W=V_{\Lambda_n}(z_n)\otimes\cdots\otimes V_{\Lambda_1}(z_1)$.

In \cite{FV2}, common eigenfunctions of $T_W(z)$, $z\in\C$
are constructed in the form
\be
b(t_1)\cdots b(t_m)v_c
\ee
where
\be
v_c(\lambda)=e^{c\lambda}\prod_{j=1}^m\frac
{\theta(\lambda-2\eta j)}
{\theta(2\eta)}
\ee
and $b(t)$ is the difference operator
$(b(t)f)(\lambda)=b(t,\lambda)f(\lambda+2\eta)$, $f\in\Fun(W)$
(both the transfer matrix and the difference operators $b(t)$
are part of the {\em operator algebra} of the elliptic
quantum group, see \cite{FV1, FV2}).
The variables $t_1,\dots,t_m$ obey a set of Bethe ansatz equations,
which are up to a shift the same as the ones described in this
paper.
The precise relation between the two approaches is the following.

\begin{thm}\label{previous}
Let $\check V=V_{\Lambda_n}(z_n)\otimes\cdots\otimes V_{\Lambda_1}(z_1)$ be
a tensor product of evaluation Verma modules with generic evaluation
points $z_1,\dots, z_n$ and let
\be
v_c(\lambda)=e^{c\lambda}\prod_{j=1}^m
\frac{\theta(\lambda-2\eta j)}{\theta(2\eta)}\in\Fun(\check V).
\ee
Then
\bea
\prod_{j=1}^mb(t_j+\eta)\,v_c
&=&
e^{c(\lambda+2\eta m)}
(-1)^m\prod_{i<j}\frac
{\theta(t_i-t_j+2\eta)}
{\theta(t_i-t_j)}\\
& &\times
\sum_{j_1+\cdots+j_n=m}\omega_{j_1,\dots,j_n}(t_1,\dots,t_n,\lambda)
\,e_{j_n}\otimes\cdots\otimes e_{j_1},
\eea
where $\omega_{j_1,\dots,j_n}(t_1,\dots,t_n,\lambda)$
is the image of $e_{j_1}^*\otimes
\cdots\otimes e_{j_n}^*$
by the map
\be
\omega(z_1,\dots,z_n,\lambda):
V_{\Lambda_1}^*\otimes\cdots\otimes V_{\Lambda_n}^*\to F_{a_1,\dots,a_n}
(z_1,\dots,z_n,\lambda).
\ee
\end{thm}
\begin{proof} The proof consists of comparing the explicit formula
of Proposition \ref{pomegam} with the explicit formula for
$\prod b(t_i)\,v_c$ given
in \cite{FV2}, Theorem 5.
\end{proof}
In particular, this gives an algebraic construction of the spaces
$F_{a}(z,\lambda)$, $a,z\in\C^n$ and their bases: we may define
$F^m_{a_1,\dots,a_n}(z_1,\dots,z_n,\lambda)$ to be the space spanned
by the coordinates of the functions
\be
\prod_{i<j}\frac
{\theta(t_i-t_j)}
{\theta(t_i-t_j+2\eta)}
\prod_{j=1}^mb(t_j+\eta)\,v_{c=0}.
\ee
of $t_1,\dots,t_m$ with values in $\check V[\,\sum\Lambda_j-2m]$.

{}From Corollary \ref{cdm} we thus obtain a completeness result for
Bethe vectors of the transfer matrix of the discrete models
of \cite{FV2}. Let $\mu$ be a generic complex number, and let
$C_\mu=\{\mu+j/N| j\in\Z\}$.
The transfer matrix $T_{\check V}(z)$ of the representation
$\check V$ is well-defined on the space $\Fun_\mu(\check V[0])$ of
functions $C_\mu\to \check V[0]$ (see \cite{FV2}, Section 4).

\begin{thm} Let $2\eta=1/N$ for some large enough odd positive integer $N$.
Let $\check V$ and $v_c$ be as in Theorem \ref{previous}, and
let $T_{\check V}(w)$ be the corresponding transfer matrix acting
on functions $C_\mu\to \check V[0]$.
\begin{enumerate}
\item[(i)] For any solution $(t_1,\dots,t_m)$ of the Bethe ansatz
equations
\be
\prod_{j:j\neq i}
\frac
{\theta(t_j-t_i-2\eta)}
{\theta(t_j-t_i+2\eta)}
\prod_{k=1}^n
\frac
{\theta(t_i-z_k-(1+\Lambda_k)\eta)}
{\theta(t_i-z_k-(1-\Lambda_k)\eta)}
=
e^{4\eta c},\qquad i=1,\dots,m,
\ee
such that, for all $i<j$, $t_i\neq t_j\mod \Z+\tau\Z$,
the vector $\psi=b(t_1)\cdots b(t_m)v\in\Fun_\mu(\check V[0])$
is a common eigenvector
of all transfer matrices $T_{\check V}(w)$ with eigenvalues
\be
\epsilon(w)=
e^{-2\eta c}\prod_{j=1}^m
\frac
{\theta(t_j-w-2\eta)}
{\theta(t_j-w)}
+
e^{2\eta c}\prod_{j=1}^m
\frac
{\theta(t_j-w+2\eta)}
{\theta(t_j-w)}
\prod_{k=1}^n
\frac
{\theta(w-z_k-(1-\Lambda_k)\eta)}
{\theta(w-z_k-(1+\Lambda_k)\eta)}
\ee
Moreover $\psi(\lambda+1)=(-1)^me^c\psi(\lambda)$.
\item[(ii)] For any generic $\alpha\in\C^\times$, there are
$d=N\dim(\check V[0])$ solutions of the Bethe ansatz equations
such that the corresponding eigenfunctions form a basis
of the space of the $d$-dimensional vector space of functions
$\psi: C_\mu\to \check V[0]$ such that $\psi(\lambda+1)=\alpha\psi(\lambda)$.
\end{enumerate}
\end{thm}
Part (i) of this Theorem is taken from \cite{FV1}. Part (ii) follows
essentially from Corollary \ref{cdm}. Note, though, that there
we had the additional hypothesis
that the $\Lambda_i$ are integers, so that the commuting operators
$H_j(z)$ are well-defined on functions on $C_\mu$.
However this hypothesis was not used
in the proof that the corresponding vectors form a basis.
Therefore the proof also applies to the present situation.

\medskip
\noi{\bf Remark.} ``Large enough'' in this theorem means
$2N>\Lambda_1+\cdots+\Lambda_n$, see Theorem \ref{toc}.
This hypothesis can probably be considerably weakened.

\medskip
\noi{\bf Remark.} The Bethe ansatz equations have ``diagonal''
solutions, i.e., solutions for which $t_i=t_j\mod \Z+\tau\Z$
for some $i\neq j$. These solutions do not in general correspond
to eigenvectors even if $b(t_1)\cdots b(t_m)v$ is finite and
nonzero. The factor
\be
\prod_{i<j}\frac
{\theta(t_i-t_j)}
{\theta(t_i-t_j+2\eta)}
\ee
in the correspondence
between the eigenvectors obtained by the algebraic Bethe ansatz
and the eigenvectors
considered in this paper sends the Bethe vectors corresponding to
diagonal solutions to zero.

\section{Solutions of the qKZB equations}\label{ssqkzb}
In this section we fix $\tau,\eta,p,\Lambda_1,\dots,\Lambda_n$, and set $a_i=\eta \Lambda_i$.
\subsection{Integral representations for solutions}

By a formal Jackson integral solution of the qKZB equations we mean
an expression
\be\Psi(z_1,\dots,z_n,\lambda)
=\int f(z_1,\dots,z_n,t_1,\dots,t_m,\lambda)Dt_1\cdots Dt_m,
\ee
where $f$ takes its values in $V[0]$,
which obeys the qKZB equations \Ref{qKZB} if we formally use the rule that
the ``integral'' $\int$ is invariant under translations of the
variables $t_i$ by $p$. In other words, $f(z_1,\dots,z_n,t_1,\dots,t_m
,\lambda)$
obeys the qKZB equations in the variables $z_i$ up to terms of the
form $g(\dots,t_i+p,\dots)-g(\dots,t_i,\dots)$.

\medskip
\begin{definition} A function $\Phi_{a}(t)$ depending on
a complex parameter $a$, such that
\be
\Phi_a(t+p)=\frac{\theta(t+a)}{\theta(t-a)}\Phi_a(t)
\ee
is called a (one-variable) phase function.
\end{definition}

We assume that $p$ has positive imaginary part, and set
$r=e^{2\pi ip}$, $q=e^{2\pi i\tau}$. Then the convergent infinite
product
\begin{equation}\label{phase1}
\Phi_a(t)=e^{-2\pi iat/p}
\prod_{j=0}^{\infty}
\prod_{k=0}^{\infty}
\frac
{(1-r^jq^ke^{2\pi i(t-a)})(1-r^{j+1}q^{k+1}e^{-2\pi i(t+a)})}
{(1-r^jq^ke^{2\pi i(t+a)})(1-r^{j+1}q^{k+1}e^{-2\pi i(t-a)})}\,,
\end{equation}
defines a phase function, and any other phase function is obtained
{}from this by multiplication by a $p$-periodic function.

Given a one-variable phase function $\Phi_a(t)$, we define
with our data an $m$-variable phase function
\begin{equation}\label{phasem}
\Phi(t_1,\dots,t_m,z_1,\dots,z_n)=
\prod_{j=1}^m
\prod_{l=1}^n\Phi_{a_l}(t_j-z_l)
\prod_{1\leq i<j\leq m}
\Phi_{-2\eta}(t_i-t_j).
\end{equation}

\begin{thm}\label{irs0} Let $\Phi_a(t)$ be a phase function,
and let $\Phi$ be the corresponding $m$-variable phase function
\Ref{phasem}. For any entire function $\xi$ of one variable, let
\be
\psi^\xi(t,z,\lambda)=\xi(p\lambda\!-\!
{\tsize\Sum_{l=1}^n2a_lz_l\!+\!4\eta\Sum_{j=1}^mt_j}))
\sum_{j_1+\cdots+j_n=m}
\omega_{j_1,\dots,j_n}(t_1,\dots,t_m,\lambda)\,
e_{j_1}\!\otimes\!\cdots\!\otimes\! e_{j_n}.
\ee
Then
\be
\Psi(z_1,\dots,z_n,\lambda)
=
\int\Phi(t_1,\dots,t_m,z_1,\dots,z_n)\psi^\xi(t_1,\dots,t_m,z_1,\dots,z_n,
\lambda)Dt_1\cdots Dt_m
\ee
is a formal Jackson integral solution of the qKZB equations.
\end{thm}

To obtain solutions from formal Jackson integral solutions, we
need to find {\em cycles}, linear forms on the space of functions
of $t_1,\dots,t_m$ that are invariant under translations $t_i\mapsto t_i+p$.
To this end we need a stronger version of the preceding theorem,
Theorem \ref{irs} below,
which gives us a space of functions on which our cycles should
be defined.

Let $\Phi$ be the phase function \Ref{phasem} and
let $a=(a_1,\dots,a_n)$, $z=(z_1,\dots,z_n)$. We assume, as
usual, that $\sum a_i=2\eta m$, $m\in\Z_{\geq 0}$.
For any entire
function $\xi$, let
$E^0_a(z;\xi)$, be the space spanned by the
functions of $t\in\C^m$ of the form
\be
\tsize \Phi(t,z)\xi(\lambda-2\eta\Sum_{k=1}^na_k+4\eta\Sum_{j=1}^mt_j)
f(t,z),
\ee
where $f(t,z)$, viewed as a function of $t=(t_1,\dots,t_m)$
belongs to $\tilde F^m_a(z,\lambda)$ (see \ref{sssym})
for some $\lambda$. All components of our integrand belong to this
space.

Let $E_a(z;\xi)$, the space of cocycles, be the space spanned by functions
of the form $g(t+p\alpha)$, where $g\in E^0_a(z;\xi)$ and $\alpha\in\Z^m$.
By construction, $E_a(z;\xi)$ is
invariant under translations of the arguments $t_i$ by $p$.
We define the space of coboundaries $DE_a(z;\xi)$
to be the subspace of $E_a(z;\xi)$
spanned by functions of the form $f(\dots, t_j+p,\dots)-
f(\dots, t_j,\dots)$, $f\in E_a(z;\xi)$.

\begin{proposition}\label{pE(z)}
$E_a(z_1,\dots,z_n)=E_a(z_1,\dots,z_j+p,\dots,z_n)$
for $j=1,\dots,n$
\end{proposition}

The proof of this proposition is part of the proof of the
following theorem, which is the main result of this section
and implies Theorem \ref{irs0}.

\begin{thm}\label{irs}
Let us write the qKZB equations as $\Psi(\dots,z_j+p,\dots)
=K_j(z)\Psi(z)$. Then, for all entire
functions $\xi$, the integrand
$\Phi(t,z)\psi^\xi(t,z)$ of Theorem \ref{irs0},
viewed as a function of $t\in\C^m$,
belongs to $E_a(z;\xi)\otimes V[0]$ for all $z\in\C^n$. It
obeys the equations
\be
\Psi(t,\dots,z_j+p,\dots)=K_j(z)\Psi(t,z)\mod DE_a(z;\xi)\otimes V[0],
\qquad j=1,\dots, m
\ee
in the cohomology $(E_a(z;\xi)/DE_a(z;\xi))\otimes V[0]$.
\end{thm}

To obtain solutions from these formal solutions, one should find horizontal
families of cycles, i.e., linear functions $\gamma(z)$ on $E_a(z;\xi)$
vanishing on $DE_a(z;\xi)$, and such that $\gamma(z+p\alpha)=\gamma(z)$
for all $\alpha\in\Z^n$. This problem
will be addressed in the next paper.

\subsection{Proof of Theorem \ref{irs}}

It is clear that the integrand belongs to $E_a(z;\xi)$.
Let us introduce some notation. If $f$ and $g$ are functions
of $t_1,\dots,t_m$ such that $\Phi f,\Phi g\in E_a(z;\xi)$, we
write $f\sim g$ if $\Phi(t,z)(f(t)-g(t))\in DE_a(z;\xi)$.
This means that $f-g$ is a linear combination of expressions
of the form $Q_ih-h$, $\Phi h\in E_a(z;\xi)$, where
\bea
&Q_if(t_1,\dots,t_m)=f(t_1,\dots,t_i+p,\dots,t_m)\phi_i(t)& \\
&\displaystyle\phi_i(t)=\prod_{l=1}^n
\frac
{\theta(t_i-z_l+a_l)}
{\theta(t_i-z_l-a_l)}
\prod_{j>i}
\frac
{\theta(t_i-t_j-2\eta)}
{\theta(t_i-t_j+2\eta)}
\prod_{j<i}
\frac
{\theta(t_i-t_j-2\eta+p)}
{\theta(t_i-t_j+2\eta+p)}
\,,\qquad i=1,\dots, m,&
\eea
The proof of the Theorem
is based on the following identity which, similarly to
Lemma \ref{l14}, follows from Lemma \ref{l13}.
\begin{lemma}
Let,
for any holomorphic function $\xi$ of one complex variable,
\be
\omega^\xi(z_1,\dots,z_n,\lambda)
=\xi(p\lambda-2\;\!(\tsize\Sum_la_lz_l-2\eta\Sum_jt_j))
\omega(z_1,\dots,z_n,\lambda).
\ee
Then, for any $\mu\in\C$ and any
$u\in V_{\Lambda_1}^*[\mu]$, $v\in (V_{\Lambda_2}^*
\otimes\cdots\otimes V_{\Lambda_{n}}^*)[-\mu]$,
\be
\prod_{j=1}^m
\frac
{\theta(t_j\!-\!z_1\!-\!a_1\!-\!p)}
{\theta(t_j\!-\!z_1\!+\!a_1\!-\!p)}
\omega^\xi(z_1+p,z_2,\dots,z_n,\lambda)u\otimes v
\sim
\omega^\xi(z_2,\dots,z_m,z_1,\lambda-2\eta\mu)u\otimes v
\ee
\end{lemma}
\begin{proof} We first assume that $\xi=1$.
We adopt the notation of Lemma \ref{l13}, and set $\tilde f=
\omega(z_1+p,\lambda)$. Then the left-hand side of
the claim is, by Lemma \ref{l13},
\begin{equation}\label{iilhs}
\frac1{(n-1)!}\Sym\biggl(\tilde f(t_1,\dots,t_{m'}) g(t_{m'+1},\dots,t_{m})
\prod_{j\leq m'} \frac{\theta(t_j-z_1-p+a_1)}{\theta(t_j-z_1-p-a_1)}\biggr),
\end{equation}
while the right-hand side is
\bea
\frac1{(n\!-\!1)!}\Sym\biggl( g(t_{m'+1},\dots,t_{m})
f(t_1,\dots,t_{m'}) \prod_{\begin{matrix}{\scriptstyle{j\leq m'}}\\
\scriptstyle{l\geq2}\end{matrix}} \frac{\theta(t_j\!-\!z_l\!+\!a_l)}
{\theta(t_j\!-\!z_l\!-\!a_l)} \prod_{\begin{matrix}{\scriptstyle{j\leq
m'}}\\ \scriptstyle{k>m'}\end{matrix}} \frac
{\theta(t_j\!-\!t_k\!-\!2\eta)} {\theta(t_j\!-\!t_k\!+\!2\eta)} \biggr).
\eea
To compare these two expressions, notice that $\tilde f=\omega(z_1+p,\lambda)u$
is related to $f=\omega(z_1,\lambda)u$ by
\be
\tilde f(t_1+p,\dots,t_{m'}+p)=f(t_1,\dots,t_{m'}).
\ee
Therefore we may rewrite the left-hand side in terms of $f$ up
to terms which are trivial in the cohomology.
For this we use the following formula, valid for general functions
$h$, which is an easy consequence
of the definition of $Q_i$:
\bea
Q_1\cdots Q_{m'}h(t_1,\dots,t_m)
&=&
h(t_1+p,\dots,t_{m'}+p,t_{m'+1},\dots,t_m)\phi^{m'}(t),\\
\phi^{m'}(t) &=&
\prod_{j=1}^{m'}
\prod_{l=1}^n
\frac
{\theta(t_j-z_l+a_l)}
{\theta(t_j-z_l-a_l)}
\prod_{i\leq m'<j}
\frac
{\theta(t_i-t_j-2\eta)}
{\theta(t_i-t_j+2\eta)}.
\eea
Therefore we have
\be
(\ref{iilhs})\sim
\frac1{(n-1)!}\Sym\biggl(f(t_1,\dots,t_{m'})
g(t_{m'+1},\dots,t_{m})
\prod_{j\leq m'} \frac{\theta(t_j-z_1+a_1)}{\theta(t_j-z_1-a_1)}
\phi^{m'}(t)\biggr),
\ee
which is easily seen to coincide with the right-hand side of
the claimed identity.

The case of general $\xi$ is proved in the same way. The point is
that the argument of $\xi$ is shifted by $-2a_1+4\eta m'$ if
$z_1$, $t_1,\dots,t_{m'}$ are shifted by $p$, while the same
effect is obtained by shifting $\lambda$ by $-2\eta\mu$.
\end{proof}

\begin{corollary}\label{cor13}
Let $1\leq j\leq n$ and let $w\in V[0]^*$ be such that
$h^{(j)}w=\mu w$. Write $z=(z_1,\dots,z_n)$, $z+p\delta_j=
(z_1,\dots,z_j+p,\dots,z_n)$, and $\rho_j(t,z)
=\prod_{i=1}^m[\theta(t_i-z_j+a_j-p)/\theta(t_i-z_j-a_j-p)]$. Then
\bea
\rho_j(t,z)\,\omega^\xi(z+p\delta_j,\lambda)
R^*_{j-1,j}(z_{j-1}\!-\!z_j\!-\!p,\lambda\!-\!2\eta\!\tsize\Sum_{l<j-1}\!
h^{(l)})\cdots\\
\cdots
R^*_{2j}(z_2\!-\!z_j\!-\!p,\lambda\!-\!2\eta h^{(1)})
R^*_{1j}(z_1\!-\!z_j\!-\!p,\lambda)w\\
\sim
\omega^\xi(z_1,\dots,z_n,\lambda\!-\!2\eta\mu)
R^*_{j,j+1}(z_j\!-\!z_{j+1},\lambda\!-\!2\eta\mu\!-\!
2\eta\tsize\Sum_{l=1}^{j-1}
h^{(l)})\cdots\\
\cdots
R^*_{j,n-1}(z_j\!-\!z_{n-1},\lambda\!-\!2\eta\tsize
\Sum_{\tsize{l=1\atop l\neq j}}^{n-2}h^{(l)})
R^*_{jn}(z_j\!-\!z_n,\lambda\!-\!2\eta\tsize
\Sum_{\tsize{l=1\atop l\neq j}}^{n-1}h^{(l)})w
\eea
\end{corollary}

This corollary follows from the above Lemma using Lemma \ref{lsym}.

The proof of Theorem \ref{irs} can now be completed: let
$w$ be a vector in $V[0]^*$ such that $h^{(j)}w=\mu w$.
Let $A(z,\lambda)=R_{1,j}(z_{1}\!-\!z_j\!-\!p,\lambda)\cdots R_{j-1,j}(z_{j-1}
\!-\!z_j\!-\!p,\lambda\!-\!2\eta\sum_{l=1}^{j-2}h^{(l)})$ and
$B(z,\lambda)=R_{jn}(z_j\!-\!z_n,\lambda\!-\!2\eta\sum_{j\neq l<n}h^{(l)})
\cdots R_{j,j+1}(z_{j}\!-\!z_{j+1},\lambda\!-\!2\eta\sum_{l<j}h^{(l)})$.
With these notations,
\bea
\langle w, A(z,\lambda)\Psi(z\!+\!p\delta_j,\lambda)\rangle
&\!=\!&
\int\Phi(t,z\!+\!p\delta_j)\,
\omega^\xi(z\!+\!p\delta_j,\lambda)A(z,\lambda)^*w
\\
&\!=\!&
\int
\Phi(t,z)\rho_j(t,z)\,
\omega^\xi(z\!+\!p\delta_j,\lambda)A(z,\lambda)^*w
\\
&\!=\!&
\int
\Phi(t,z)\,
\omega^\xi(z,\lambda\!-\!2\eta\mu)B(z,\lambda\!-\!2\eta\mu)^*w
\\
&\!=\!&\langle w, B(z,\lambda\!-\!2\eta\mu)\Psi(z,\lambda\!-\!2\eta\mu)\rangle.
\eea
The first equality follows from the definition of $\Psi$. The
second is the identity $\Phi(t,z+p\delta_j)=\rho_j(t,z)\Phi(t,z)$,
an immediate consequence of the defining property of the phase
function. The third equality follows from Corollary \ref{cor13}
and the last equality by the definition of $\Psi$.
Thus $A(z,\lambda)\Psi(z+p\delta_j,\lambda)
=\Gamma_jB(z,\lambda)\Psi(z,\lambda)$.
To complete the proof of the Theorem, one has to invert the operator $A$,
which is very easy, given the unitarity of $R$-matrices
(property III of Section \ref{sqkzb}).

\section{Residues}\label{sp}
In this section we give a proof of Proposition \ref{p32}, and
some other technical results. For this we introduce the
basic computational tool of (iterated) residues.

Let us first describe the spaces
$F_{a_1,\dots,a_n}^m(z_1,\dots,z_n,\lambda)$ in terms of symmetric
theta functions. Let $G^m_{a_1,\dots,a_n}(z_1,\dots,z_n,\lambda)$ be
the space of meromorphic functions $g$ of $m$ complex
variables $t_1,\dots, t_m$ such that
\begin{enumerate}
\item[(i)] $\prod_{i=1}^m\prod_{k=1}^n\theta(t_i\!-\!z_k\!-\!a_k)g
(t_1,\dots,t_m)$
is a holomorphic function of $t_1,\dots,t_m\in\C$.
\item[(ii)] $g$ is periodic with period 1 in each of its arguments
and
\be
g(\dots,t_j\!+\!\tau,\dots)=e^{-2\pi i(\lambda
\!+\!2\eta m)}g(\dots,t_j,\dots).
\ee
for all $j=1,\dots,m$.
\item[(iii)] $g(t_{\sigma(1)},\dots,t_{\sigma(m)})=
g(t_1,\dots,t_m)$ for any permutation $\sigma\in S_m$.
\end{enumerate}

\begin{lemma}\label{limo}
The linear map sending $f$ to
\be
g(t_1,\dots,t_m)=\prod_{i<j}\frac{\theta(t_i\!-\!t_j\!+\!2\eta)}
{\theta(t_i\!-\!t_j)}f(t_1,\dots,t_m)
\ee
is an isomorphism from
$F^m_{a_1,\dots,a_n}(z_1,\dots,z_n,\lambda)$
to
$G^m_{a_1,\dots,a_n}(z_1,\dots,z_n,\lambda)$
\end{lemma}

\begin{proof} By using the transformation properties of
$\theta$ under translation by $\Z+\tau\Z$, and the fact that
$\theta$ is odd, we see that $f$ has the required behavior
under lattice translation and the action of permutations if
and only if $g$ does.
Since the poles of $\theta$ lie in the lattice, it is
also clear that if $g$ obeys (i), then the poles
of $f$ are on the correct hyperplanes.
Conversely, if $f\in F^m_{a_1,\dots,a_n}(z_1,\dots,z_n,\lambda)$,
then $\prod_{i<j}\theta(t_i\!-\!t_j\!+\!2\eta)f(t_1,\dots,t_m)$ is regular
except possibly on the hyperplanes $t_i\equiv z_j\!+\!a_j$ (modulo
the lattice). From the symmetry properties of $f$ and
the fact that $\theta$ is odd, it follows that this function
is skew-symmetric under permutation of the
arguments. In particular, it vanishes on the diagonals $t_i\!-\!t_j=0$,
and, by the transformation property of $f$ under translation by $\Z+\tau\Z$,
also on the hyperplanes $t_i\!-\!t_j=r\!+\!s\tau$, $r$, $s\in\Z$.
Therefore we can divide by $\prod_{i<j}\theta(t_i\!-\!t_j)$ without
creating new poles. Thus $g$ obeys (i).
\end{proof}

\begin{lemma}\label{ldim}
For any $a_1,\dots,a_n,z_1,\dots,z_n,\lambda\in\C$ and $m\in\Z_{\geq 0}$
the dimension of the space
$F_{a_1,\dots,a_n}(z_1,\dots,z_n,\lambda)$
is
\be
\left(\begin{matrix}{n+m-1}\\ {m}\end{matrix}\right),
\ee
the number of ways of decomposing $m$ as
a sum of $n$ nonnegative integers.
\end{lemma}
\begin{proof}
Multiplication by $\prod_{i=1}^m\prod_{j=1}^n\theta(t_i\!-\!z_j\!-\!a_j)$
is an isomorphism from the space
$G^m_{a_1,\dots,a_n}(z_1,\dots,z_n,\lambda)$
to the space of entire functions that are invariant under
permutations of the arguments and have theta function transformation
properties under lattice translation of each argument. The
dimension of the latter space is easily computed, say by Fourier
series.
\end{proof}

In particular, if $n=1$, we have:

\begin{lemma}\label{ln1}
$F^m_a(z,\lambda)$ is a one-dimensional space spanned by
\be
\omega_m(t_1,\dots,t_m)=
\prod_{i<j}\frac{\theta(t_i\!-\!t_j)}
{\theta(t_i\!-\!t_j\!+\!2\eta)}\prod_{j=1}^{m}\frac
{\theta(\lambda\!+\!2\eta m\!+\!t_j\!-\!z\!-\!a)}
{\theta(t_j\!-\!z\!-\!a)},
\ee
\end{lemma}

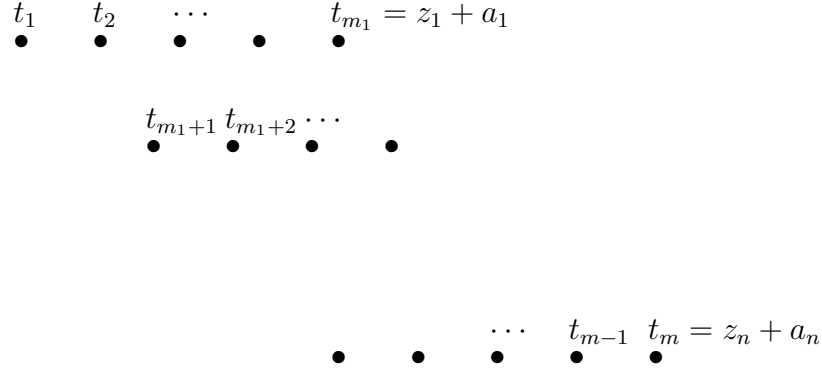
\begin{figure}
\begin{picture}(350,200)
\put(50,150){\makebox{$\bullet$}}
\put(50,160){\makebox{$t_1$}}
\put(80,150){\makebox{$\bullet$}}
\put(80,160){\makebox{$t_2$}}
\put(110,150){\makebox{$\bullet$}}
\put(110,160){\makebox{$\cdots$}}
\put(140,150){\makebox{$\bullet$}}
\put(170,150){\makebox{$\bullet$}}
\put(170,160){\makebox{$t_{m_1}=z_1+a_1$}}
\put(100,110){\makebox{$\bullet$}}
\put(100,120){\makebox{$t_{m_1+1}$}}
\put(130,110){\makebox{$\bullet$}}
\put(130,120){\makebox{$t_{m_1+2}$}}
\put(160,110){\makebox{$\bullet$}}
\put(160,120){\makebox{$\cdots$}}
\put(190,110){\makebox{$\bullet$}}
\put(170,30){\makebox{$\bullet$}}
\put(200,30){\makebox{$\bullet$}}
\put(230,40){\makebox{$\cdots$}}
\put(230,30){\makebox{$\bullet$}}
\put(260,40){\makebox{$t_{m-1}$}}
\put(260,30){\makebox{$\bullet$}}
\put(290,40){\makebox{$t_m=z_n+a_n$}}
\put(290,30){\makebox{$\bullet$}}
\end{picture}
\caption{The point at which the residue is taken.}
\label{fig1}
\end{figure}

The next tool is a suitable set of iterated residues spanning, for
generic parameters, the vector
space dual to $F^m_{a_1,\dots,a_n}(z_1,\dots,z_n,\lambda)$.
For any set $m_1,\dots,m_n$ of nonnegative integers with
sum $m$ and
any meromorphic function $f(t_1,\dots,t_m)$ in $m$ variables
we let $\res_{m_1,\dots_{m_n}}f$ be the complex number
\bea
\res_{m_1,\dots,m_n}f=
\res_{t_1=z_1+a_1-2\eta(m_1-1)}
\cdots
\res_{t_{m_1-1}=z_1+a_1-2\eta}
\res_{t_{m_1}=z_1+a_1}
\\
\res_{t_{m_1+1}=z_2+a_2-2\eta(m_2-1)}
\cdots
\res_{t_{m_1+m_2-1}=z_2+a_2-2\eta}
\res_{t_{m_1+m_2}=z_2+a_2}\\
\cdots
\res_{t_{m-1}=z_n+a_n-2\eta}
\res_{t_{m}=z_n+a_n}f.
\eea
The map $f\to \res_{m_1,\dots,m_n}f$ is a linear function
on $F^m_{a_1,\dots,a_n}(z_1,\dots,z_n,\lambda)$. The point
at which this residue is taken is represented in Fig.\ \ref{fig1}

\begin{proposition}\label{pres}
Fix $a_1,\dots,a_n,z_1,\dots,z_n,\lambda$.
\begin{enumerate}
\item[(i)] If $z_1,\dots,z_n,\lambda$ and $\eta$ are generic, then
the residues $\res_{m_1,\dots,m_n}$
with $m_1+\cdots+m_n=m$ form a basis of the dual of
$F^m_{a_1,\dots,a_n}(z_1,\dots,z_n,\lambda)$.
\item[(ii)] Let $\omega_{j_1,\dots,j_n}\in
F^m_{a_1,\dots,a_n}(z_1,\dots,z_n,\lambda)$
denotes the image by $\Phi_n$ of
the basis $\omega_{j_1}\otimes\cdots\otimes\omega_{j_n}$
of $\otimes_{l=1}^n
F^{j_l}_{a_l}(z_l,\lambda-2\sum_{k=1}^{l-1}(a_k-2\eta j_k))$,
then
\be
\res_{m_1,\dots,m_n}\omega_{j_1,\dots,j_n}=0
\ee
unless $m_l+\cdots+m_n\leq j_l+\cdots+j_n$ for all $l=1,\dots,n$.
\item[(iii)]
$\res_{m_1,\dots,m_n}\omega_{m_1,\dots,m_n}\neq 0$
for generic $z_1,\dots,z_n,\lambda$ and $\eta$.
More precisely,
\bea
\res_{m_1,\dots,m_n}\omega_{m_1,\dots,m_n}
&\!=\!&\prod_{i=1}^n\prod_{j=1}^{m_i}
\frac{\theta(\lambda-2\sum_{l<i}(a_l-2\eta m_l)+2\eta j)}
{\theta'(0)}
\\
& &
\prod_{1\leq l<k\leq n}\prod_{j=1}^{m_k}
\frac
{\theta(z_k+a_k-z_l+a_l-2\eta(j-1))}
{\theta(z_k+a_k-z_l-a_l-2\eta(j-1))}
\eea
\end{enumerate}
\end{proposition}

\begin{proof}
We prove (ii) and (iii), from which (i) follows.

Let first $n=1$. By Lemma \ref{ldim}, $F_a(z,\lambda)$ is
one dimensional. We compute:
\bea
\res_{m}\omega_m&\!=\!&
\theta'(0)^{-m}\frac{\prod_{i<j}{\theta(2\eta(i-j))}}
{\prod_{i<j-1}{\theta(2\eta(i-j)+2\eta)}}
\frac{\prod_{j=1}^m{\theta(\lambda+2\eta m-2\eta(m-j))}}
{\prod_{j=1}^{m-1}\theta(-2\eta(m-j))}
\\
&\!=\!&
\prod_{j=1}^{m}\frac{\theta(\lambda+2\eta j)}{\theta'(0)}.
\eea
This does not vanish for generic $\lambda$.

Consider now the general case.

We introduce the convenient notation
$t_I=(t_{i_1},\dots,t_{i_m})$ for any finite set of integers
$I=\{i_1<\cdots<i_m)$. Then $\omega_{j_1,\dots,j_n}$
can be written as
\be
\sum_{I_1,\dots, I_n}\omega_{j_1}(t_{I_1})\cdots
\omega_{j_n}(t_{I_n})
\prod_{l<k}\prod_{i\in I_k}
\frac{\theta(t_i-z_l+a_l)}
{\theta(t_i-z_l-a_l)}
\prod_{i<j,k>l}
\prod_{i\in I_k,j\in I_l}
\frac{\theta(t_i-t_j-2\eta)}
{\theta(t_i-t_j+2\eta)}.
\ee
The sum is taken over all partitions $I_1,\dots,I_n$ of
$\{1,\dots,m\}$ such that $I_l$ contains precisely $j_l$
elements, $1\leq l\leq m$.

The next step is to determine the partitions corresponding
to terms with non-trivial residue $\res_{m_1,\dots,m_n}$.
We may assume that all parameters are generic. Then
a term with non-trivial residue must in particular
have a pole on the hyperplane $t_m=z_n+a_n$, if $m_n>0$.
This implies that $m\in I_n$. If $m_n>1$, then,
we must also have a pole on the hyperplane $t_{m-1}-t_m+2\eta=0$.
But a pole on the hyperplane $t_j-t_{j+1}+2\eta$ occurs
only if $j$, $j+1$ belong to the same $I_k$, (the pole
is in $\omega_k$) or if $j\in I_k$, $j+1\in I_l$ for
some $l<k$ (the pole is in the last factor in this case).
Thus $m-1\in I_m$, and so on. We obtain the necessary condition
$\{m_1+\cdots+m_{n-1}+1,\dots,m\}\subset I_n$. The next residue
is at $t_{m_1+\cdots+m_{n-1}}=z_{n-1}+a_{n-1}$. Thus
we see that $m_1+\cdots+m_{n-1}$ must belong to $I_{n-1}\cup I_{n}$.
Continuing
this way, we conclude
that a partition $I_1,\dots,I_n$ corresponds to a term
with trivial residue unless
\be
\{m_1+\cdots+m_l+1,\dots,m\}\subset I_{l+1}\cup\cdots\cup I_n,
\qquad
l=1,\dots,n-1.
\ee
In particular this proves (ii).

If $m_l=j_l$ for all $l$, then only one partition contributes
to the residue, namely $I_1=\{1,\dots,m_1\}$, $I_2=\{m_1+1,
\dots,m_1+m_2\},\dots$ The residue is then the product of
the residues of the individual $\omega_{j_i}$'s times the
value of the remaining product
product of ratios of theta functions, which can easily
be computed.

The claim (i) follows from (ii), (iii) and Lemma \ref{ldim}.
\end{proof}

\noi{\it Proof of Proposition \ref{p32}}:
By Proposition \ref{pres}, $\Phi_n$ is generically an isomorphism.
{}From this and the associativity of $\Phi$, the main
claim of Proposition \ref{p32} follows.
$\;\Box$
\vs

\noi{\it Proof of Lemma \ref{lpp}}, (ii):
We have to show that $\omega(w,z,\lambda):
V_\Lambda ^*\otimes V_\Lambda ^*\to F_{aa}(z,w,\lambda)$ is
invertible at $z=w$. Since the dimension of
$F_{aa}(z,w,\lambda)$
is $m+1$
for all $z,w$, it is sufficient to show that, for
any generic $a$, the functions $\omega_{j,m-j}=\omega(z,z,\lambda)e_j^*\otimes
e_{m-j}^*$, $0\leq j\leq m$, are linearly independent.
Here we have double poles, the residues $\res_{j,m-j}$ are no
longer linearly independent and we have to change
slightly the construction. Let
$\tilde\res_{m,0}=\res_{m}$ and if $j=1,\dots,m$,
let $\tilde\res_{j,m-j}f$ be
\bea
&
\res_{t_1=z+a-2\eta(j-1)}
\cdots
\res_{t_{j-1}=z+a-2\eta}\res_{t_j=z+a}
&
\\
&
\res_{t_{j+1}=z+a-2\eta(m-j-1)}
\cdots
\res_{t_{m-1}=z+a-2\eta}
\res_{t_{m}=z+a}[(t_m-z-a)f].
&
\eea
Then, by proceeding the same way as in the proof of Proposition \ref{pres}
we see that the matrix
$(\tilde\res_{j,m-j}\omega_{k,m-k})_{j,k}$
is triangular with non-vanishing diagonal
matrix elements.
Therefore $\omega(z,z,\lambda)$ is
invertible, and
\be\lim_{w\to z}R^*_{\Lambda\Lambda}(z,w,\lambda)=
\omega(z,z,\lambda)\omega(z,z,\lambda)^{-1}P=P.
\ee
$\;\Box$

\subsection{Resonances}

Here we consider the case when the parameters
$\Lambda_i$ are nonnegative integers.
{}From the point of view of elliptic
quantum groups (see \cite{FV1}), this is the case when
the evaluation Verma module has a finite dimensional quotient.
In the present framework, what happens is that the
$R$-matrix has remarkable invariant subspaces.

Let us introduce some notation: if $\ell\in\Z_{\geq 0}$,
we let $L^*_\ell=\oplus_{j=0}^\ell\C e_j^*$, a finite dimensional
subspace of $V^*_\ell$. Thus its dual space
$L_\ell$ is a quotient of $V_\ell$, which
we may identify with the subspace $\oplus_{j=0}^\ell \C e_j$.

If $f$ is a meromorphic function of $m$ variables
$t_1,\dots,t_m$, and $z\in\C$ we define
for every integer $\ell\in\{0,\dots,m-1\}$ a function $\res_{(z,\ell)}f$
of $t_1,\dots,t_{m-\ell-1}$:
\be
\res_{(z,\ell)}f=\res_{t_{m-\ell}=z-2\eta \ell}\cdots
\res_{t_{m-1}=z-2\eta}\res_{t_{m}=z}.
\ee
We extend this definition to all $\ell\in\Z_{\geq 0}$ by
setting $\res_{(z,\ell)}=0$ if $\ell\geq m$.

Now consider the map $\omega(z,w,\lambda): V^*_\Lambda \otimes V^*_\Mu\to
F_{ab}(z,w,\lambda)$ defined in Section \ref{ssafs}. We
let $\Ker\,\res_{(u,\ell)}$ be the space of $f\in F_{ab}(z,w,\lambda)$
such that $\res_{(u,\ell)}f=0$.

\begin{thm}\label{treso}
Let $z,w,\eta,\lambda$ be generic and $a=\eta\Lambda$, $b=\eta\Mu\in\C$.
\begin{enumerate}
\item[(i)] If $\Lambda\in\Z_{\geq 0}$, then
$\omega(z,w,\lambda)$ maps $L_\Lambda ^*\otimes V_\Mu^*$ onto
$\Ker\,\res_{(z+a,\Lambda)}$.
\item[(ii)] If $\Mu\in\Z_{\geq 0}$, then
$\omega(z,w,\lambda)$ maps $V_\Lambda ^*\otimes L_\Mu^*$ onto
$\Ker\,\res_{(w+b,\Mu)}$.
\item[(iii)] If $\Lambda\in\Z_{\geq 0}$ and $\Mu\in
\Z_{\geq 0}$, then
$\omega(z,w,\lambda)$ maps $L_\Lambda ^*\otimes L_\Mu^*$ onto
$\Ker\,\res_{(z+a,\Lambda)}\cap\Ker\,\res_{(w+b,\Mu)}$.
\end{enumerate}
\end{thm}

\begin{proof}
We prove (i). The other claims are proved similarly, or in
a simpler way.

Thus we assume that $\Lambda$ is a nonnegative integer $\ell$.
We first show that the image of $L^*_\Lambda \otimes V^*_\Mu$ is
contained in the kernel of the residue. In other words,
we must show that if $f\in F^{m'}_a (z,\lambda)$ and
$g\in F^{m''}_b(w,\lambda-2a+4\eta m')$, with $m'\leq \ell$,
then $\res_{(z+a,\ell)}\Phi(f\otimes g)=0$. Set $m'+m''=m$.
The function $\Phi(f\otimes g)$ is, in the notation used in the
proof of Proposition \ref{pres}
\be
\sum_{I_1,I_2}f(t_{I_1})g(t_{I_2})
\prod_{i\in I_2}
\frac{\theta(t_i-z+a)}
{\theta(t_i-z-a)}
\prod_{i\in I_2,j\in I_1,i<j}
\frac{\theta(t_i-t_j-2\eta)}
{\theta(t_i-t_j+2\eta)}.
\ee
The sum runs over all partitions of $\{1,\dots,m\}$ into
pairs of disjoint subsets $I_1$, $I_2$ of cardinality $m'$,
$m''$, respectively.

The partitions with non-vanishing residue $\res_{(z,\ell)}$
must correspond to terms with poles on the hyperplanes
$t_{m-\ell}-t_{m-\ell+1}+2\eta=0, \cdots, t_{m-1}-t_m+2\eta=0$. Thus
if $m-1\geq j\geq m-\ell$ and $j\in I_k$ $(k=1,2)$, then $j+1\in I_r$ with
$r\leq k$. Since $m'\leq \ell$, this implies that $m-\ell\in I_2$.
But then the last residue $\res_{t_{m-\ell}=z+a-2\eta\ell}$
vanishes since the pole at $z+a-2\eta\ell=z-a$
is canceled by the zero of
$\theta(t_{m-\ell}-z+a)$.

To show that the image is precisely the kernel of the residue, we
proceed as in the proof of
Proposition \ref{pres}: suppose $f\in F^m_{ab}(z,w)$ is in
$\Ker\;\res_{(z+a,\ell)}$. We claim that $f$ is in the image of
$L^*_\Lambda \otimes V^*_\Mu$. Let us think of $f$ as a linear combination of
$\omega_{j,m-j}$, $j\in\{0,\dots,m\}$. Consider the residues
\be
\widehat\res_{m-k,k}=\res_{t_1=w+a-2\eta(m-k-1)}\cdots\res_{t_{m-k}=w+a}
\res_{t_{m-k+1}=z+a-2\eta(k-1)}\cdots\res_{t_m=z+a}.
\ee
As in the proof of Proposition \ref{pres}, we see that
$\widehat\res_{m-k,k}\omega_{j,m-j}$ vanishes if $k<j$, and is
non-zero for generic $z,w,\lambda$, if $k=j$. Also if
$f\in\Ker\,\res_{(z+a,\ell)}$, then $\widehat\res_{m-k,k}f=0$
whenever $k>\ell$. In particular, $\widehat\res_{m-k,k}\omega_{j,m-j}=0$
if $j\leq\ell$ and $k>\ell$. Therefore, $\widehat\res_{m-\ell-1,\ell+1}f=0$
implies that the coefficient of $\omega_{\ell+1,m-\ell-1}$ vanishes;
this and $\widehat\res_{m-\ell-2,\ell+2}f=0$ implies that the coefficient
of $\omega_{\ell+2,m-\ell-2}$ vanishes, and so on. We
conclude that $f$ is a linear combination of $\omega_{j,m-j}$
with $j\leq \ell$, which proves our claim.
\end{proof}

\end{document}